\documentstyle[twoside,12pt]{article}
\newtheorem{prop}{Proposition}[section]
\newtheorem{dfn}[prop]{Definition}
\newtheorem{theo}[prop]{Theorem}
\newtheorem{conj}[prop]{Conjecture}
\newtheorem{rem}[prop]{Remark}
\newtheorem{coro}[prop]{Corollary}

\newtheorem{exam}[prop]{Example}
\newtheorem{ques}[prop]{Question}

\title{\bf Strong McKay Correspondence, \\
String-theoretic  Hodge Numbers \\
and Mirror Symmetry}

\author{ Victor V. Batyrev\thanks{ Supported by DFG.} \\
\small FB Mathematik,  Universit\"at-GHS-Essen \\
\small  Universit\"atsstr\ss e 3,  45141  Essen, FRG  \\
\small e-mail: victor.batyrev@aixrs1.hrz.uni-essen.de \\
 and \\
Dimitrios I. Dais \\
\small Max-Plank-Institut f\"ur Mathematik \\
\small  Gottfried-Claren-Str. 26, 53225 Bonn, FRG  \\
\small e-mail: dais@mpim-bonn.mpg.de}

\begin{document}

\date{}

\maketitle

\thispagestyle{empty}

\begin{abstract}
We propose a new conjectural version of the  McKay
correspondence which enables us to understand
the ``Hodge numbers'' assigned to singular
Gorenstein varieties by physicists.
Our results lead to the conjecture that string theory indicates
the existence of some new cohomology theory $H^*_{\rm st}(X)$
for algebraic varieties with  Gorenstein singularities. We give
a formal mathematical definition of the Hodge numbers  $h^{p,q}_{\rm st}(X)$
inspired from the consideration of strings on orbifolds and
from this new conjectural version of the McKay correspondence.  The
numbers $h^{p,q}_{\rm st}(X)$ are expected to give
the spectrum of orbifoldized
Landau-Ginzburg models and mirror duality relations
for  higher dimensional Calabi-Yau varieties with Gorenstein toroidal or
quotient singularities.
\end{abstract}

\newpage

\section{Introduction}

Throughout this paper by an {\em algebraic variety} (or simply {\em
variety})  we mean an integral, separated algebraic scheme over ${\bf C}$.
By a {\em compact algebraic variety} we mean the representative of a
complete variety within the analytic category. The
{\em singular} {\em locus} of an algebraic variety
$X$ is denoted by ${\rm Sing}\,X$. The
words {\em smooth variety} and {\em manifold} are used interchangeably. By
the word {\em singularity} we sometimes intimate a singular point and
sometimes
the underlying space of a neighbourhood or the germ of a singular point,
but its meaning will be always clear from the context. Following Danilov
\cite{danilov}, \S 13.3, we shall say that an $x \in X$ is a
{\em toroidal singularity} of $X$, if there is an analytic isomorphism
between the germ $(X, x)$ and the germ corresponding to the toric singularity
$({\bf A}_{\sigma}, p_{\sigma})$ (see also \S 4).

Our main tool will be certain algebraic varieties
with special Gorenstein singularities, primarily having in mind
the Calabi-Yau varieties. A {\em Calabi-Yau
variety} is defined to be a normal projective algebraic variety $X$
with trivial canonical sheaf ${\omega}_X$ and $H^i(X, {\cal O}_X) = 0$,
$0 < i < {\rm dim}_{\bf C}\,X$, which, in addition, can have at most
{\em canonical Gorenstein singularities}. (For the notion of
{\em canonical singularity} we refer to \cite{reid1}.)  If ${\rm Sing}\,X =
\emptyset$, then $X$ is called, as usual, {\em Calabi-Yau manifold}.

In this paper we shall attempt to realize  some Hodge-theoretical
invariants used by physicists for
singular varieties being related  to the mirror symmetry phenomenon.
The necessity of working with singular varieties becomes unavoidable
from the fact that, in many examples of pairs $X$, $X^*$ of mirror symmetric
Calabi-Yau manifolds, at least one of the two manifolds $X$ or $X^*$ is
obtained as {\em a crepant desingularization} of a singular
Calabi-Yau variety \cite{batyrev1,morrison}.
Here, by a crepant desingularization of a
Gorenstein variety $Z$, we mean a birational morphism $\pi\,: \,
Z' \rightarrow Z$, such that $\pi^*(\omega_Z) \cong \omega_{Z'}$,
where $\omega_Z$ and $\omega_{Z'}$ denote the canonical sheaves
on $Z$ and $Z'$ respectively. 3-dimensional Gorenstein
quotient singularities and their crepant desingularizations
have been studied in
\cite{bertin,ito1,ito2,markushevich,markushevich2,roan0,roan1,roan2,roan3,roan-y
au,reid2,yau}.

The most known physical cohomological invariant of
singular varieties  obtained as
quotient-spaces of certain compact manifolds by actions of finite groups
is the so called
{\em physicists Euler number} \cite{dixon}. It has been investigated
by several mathematicians in  \cite{atiyah,got2,hirzebruch,roan,roan1,reid2}.

Let $X$ be a smooth simplectic manifold over ${\bf C}$ having
an action of a finite group $G$ such that the simplectic
volume form $\omega$ is $G$-invariant. For any $g \in G$, we set $X^g
: = \{ x \in X \mid g(x) = x \}$.
Physicists have proposed the following  formula  for computing
the {\em orbifold Euler number}  \cite{dixon}:
\begin{equation}
 e(X,G) = \frac{1}{\mid G \mid} \sum_{gh = hg} e(X^g \cap X^h).
\label{euler.phys}
\end{equation}
It is expected that $e(X,G)$ coincides with the usual Euler number
$e(\widehat{X/G})$ of a crepant desingularization $\widehat{X/G}$ of
the quotient space $X/G$ provided  such a desingularization exists.
For a volume-invariant  linear action on ${\bf C}^n$ of a finite group $G$,
the corresponding  conjectural local  properties
of crepant desingularizations were formulated by M. Reid \cite{reid2}:

\begin{conj}
{\rm (generalized McKay correspondence)} Let $X = {\bf C}^n$,
$G$ an arbitrary finite subgroup in $SL(n, {\bf C})$. Assume that
$Y = X/G$ admits a crepant desingularization $\pi \,: \, \hat{Y}
\rightarrow Y$. Then $H^*(\pi^{-1}(0), {\bf C})$ has a basis consisting of
classes of algebraic cycles $Z_c \subset \pi^{-1}(0)$ which are
in $1$-to-$1$ correspondence with conjugacy classes $c$ of $G$.
In particular, we obtain for the Euler number
\[ e(\hat{Y}) = e(\pi^{-1}(0))
 = \# \{ \mbox{\rm conjugacy classes in $G$} \}. \]
 \label{general}
\end{conj}

\begin{rem}
{\rm For $n =2$ an one-to-one   correspondence between the nontrivial
irreducible representations of a subgroup $G \subset SL(2, {\bf C})$ and
the irreducible components of $\pi^{-1}(0)$ was discovered by McKay
\cite{mckay} and investigated in \cite{gonzalez,knorrer,sandro-infirri}. }
\end{rem}

Our first purpose is to use some stronger version of Conjecture \ref{general}
in order to give an analogous
interpretation for the {\em physicists Hodge numbers} $h^{p,q}(X,G)$
of orbifolds considered by C. Vafa \cite{vafa} and E. Zaslow \cite{zaslow}.
Let $X$ be a smooth compact
K\"ahler  manifold of dimension $n$ over ${\bf C}$ being equipped with
an action of a finite group $G$,  such that $X$ has a $G$-invariant
volume form. Let $C(g) : = \{ h \in G \mid hg = gh\}$. Then the action
of $C(g)$ on $X$ can be restricted on $X^g$. For any point $x \in X^g$, the
eigenvalues of $g$ in the holomorphic tangent space $T_x$ are roots of unity:
\[ e^{2\pi \sqrt{-1} \alpha_1}, \ldots, e^{2\pi \sqrt{-1} \alpha_d} \]
where $0 \leq \alpha_j <1 $ $(j =1, \ldots, d)$ are locally constant
functions on $X^g$ with values in ${\bf Q}$.
We write  $X^g = X_1(g) \cup \cdots \cup
X_{r_g}(g)$, where $X_1(g), \ldots, X_{r_g}(g)$ are the smooth  connected
components of $X^g$. For each $ i \in \{1, \ldots, r_{g} \}$,
the {\em fermion shift number} $F_i(g)$ is defined
to be equal to the value of  $\sum_{1 \leq j \leq n}
{\alpha_j}$ on the connected component $X_i(g)$.
We denote by $h^{p,q}_{C(g)}(X_i(g))$ the dimension of the subspace of
$C(g)$-invariant elements in $H^{p,q}(X_i(g))$. We set
\[ h^{p,q}_g(X,G) : = \sum_{i =1}^{r_g} h^{p - F_i(g),q -
F_i(g)}_{C(g)}(X_i(g)).\]
The  {\em orbifold Hodge numbers} of $X/G$ are defined by the
formula (3.21) in \cite{zaslow}:
\begin{equation}
 h^{p,q}(X,G): = \sum_{\{ g \}} h^{p,q}_g(X,G)
 \label{phys.form}
\end{equation}
where $\{g\}$ runs over  the conjugacy classes of $G$, so that $g$ represents
$\{ g \}$. As we shall see in Corollary \ref{c.des},
these numbers coincide with
the {\em usual} Hodge numbers of a crepant desingularization of $X/G$.

One of our next intentions is to convince the reader of the existence
of some {\em new cohomology theory} $H^*_{\rm st}(X)$ of more general
algebraic varieties $X$ with mild Gorenstein singularities.
Since this cohomology is inspired from the  string
theory,  we call $H^*_{\rm st}(X)$ the {\em string cohomology of} $X$.
For  compact varieties $X$, we expect
 that the string cohomology groups $H^*_{\rm st}(X)$ will satisfy the
 Poincar\'{e} duality and will be endowed with a pure Hodge structure.
The role  of   crepant resolutions for   the string-cohomology
$H^*_{\rm st}(X)$ is analogous to that one  of small resolutions
for the intersection cohomology $IH^*(X)$  with  middle perversity.
Physicists  compute orbifold Hodge numbers without using
crepant  desingularizations. From mathematical point of view, however,
crepant  desingularizations seem to  be rather helpful, although  they
have some  disadvantages. Firstly, they might not exist
(at least in dimension $\geq 4$) and ,secondly, even if they
exist, they might be not unique.
The consistency of the physical approach naturally suggests
the formulation of the following
conjecture (which  can be verified for the toric case by
Theorem \ref{invariants}):

\begin{conj}
Hodge numbers of smooth crepant resolutions do not depend on the choice of
such a resolution.
\end{conj}

Let us briefly review the rest of the paper. In Section 2, we consider an
example showing the importance  of
 the ``physical Hodge numbers''  in connection with
 the mirror duality.
 In Section 3,  we remind basic properties of
$E$-polynomials. In Section 4,  we study the Hodge structure of
the exceptional loci of local crepant toric resolutions. In Section 5,
we formulate
the conjecture concerning the strong McKay correspondence and we prove that
it is true for $2$- and $3$-dimensional Gorenstein quotient singularities,
as well as for abelian Gorenstein quotient singularities of arbitrary
dimension.
This  correspondence will be used in Section
6 in order to give the  formal definition of
the string-theoretic Hodge numbers and to study their main properties.
In Section
7, we give some applications relating to the mirror symmetry and
formulate the string-theoretic Hodge diamond-mirror conjecture
for  Calabi-Yau complete intersections in $d$-dimensional
toric Fano varieties.  This conjecture will be proved
in Section 8 for the case of $\Delta$-regular hypersurfaces in
toric Fano varieties ${\bf P}_{\Delta}$ which are defined by
$d$-dimensional reflexive simplices $\Delta$ (for arbitrary $d$); it gives
the mirror duality {\em for all string-theoretic Hodge numbers}
$h^{p,q}_{\rm st}$ of abelian quotients of
Calabi-Yau Fermat-type hypersurfaces which are embedded in $d$-dimensional
weighted projective spaces.  This  duality agrees  with  the mirror
construction proposed by Greene and Plesser \cite{greene0,greene1,greene}
and the polar duality of reflexive polyhedra proposed in \cite{batyrev1}.
\bigskip

{\bf Acknowledgements.} We would like to express our thanks to D. Cox, A.
Dimca,
H. Esnault, L. G\"ottsche, Yu. Ito, D. Kazhdan, M. Kontsevich,
D. Markushevich,
Yu. Manin, K. Oguiso,
M. Reid, A. V. Sardo-Infirri, D. van Straten and E. Viehweg
for fruitful discussions, suggestions and remarks.

\section{Hodge numbers
and mirror symmetry}

At the beginning we shall state some introductory questions
which could be considered also
as another  motivation for the paper. These questions are related
to singular varieties  of dimension $\geq 4$ which arose
as examples of the mirror duality
\cite{batyrev1,candelas,greene,schimmrigk1,schimmrigk2,schimmrigk3}.
If two $d$-dimensional Calabi-Yau manifolds $X$ and $Y$
form a mirror pair, then  for all $0 \leq p, q \leq d$
their  Hodge numbers must satisfy
the relation
\begin{equation}
 h^{p,q}(X) = h^{d-p,q}(Y).
\label{duality}
\end{equation}
However, it might happen that a mirror pair consists of two
$d$-dimensional Calabi-Yau varieties $X$ and $Y$ having singularities. In
this case, the duality (\ref{duality}) is expected to take place
not for $X$ and $Y$
themselves, but for their crepant desingularizations $\hat{X}$ and
$\hat{Y}$, if such desingularizations exist. Using the existence of smooth
crepant desigularizations of Gorenstein toroidal singularities
in dimension $\leq 3$, one can check the relations (\ref{duality})
for many examples of $3$-dimensional mirror pairs
\cite{batyrev1,roan0}. But there are  difficulties to prove
(\ref{duality}) {\em for all} $p,q$ and  $d \geq 4$, even if
one heuristically
knows a mirror pair of singular Calabi-Yau varieties, for instance, as
an orbifold.  The main problem in dimension
$d \geq 4$ is due to the existence  of many  {\em terminal} Gorenstein
quotient singularities, i.e., to singularities
which obviously do not admit any crepant resolution.
In \cite{batyrev1}, the first author constructed the  so called {\em maximal
projective crepant partial desingularizations} (MPCP-desingularizations)
of singular Calabi-Yau hypersurfaces in toric varieties.
Using MPCP-desingularizations,  the relation
(\ref{duality}) was proved  for $h^{1,1}$ and $h^{d-1,1}$ in \cite{batyrev1}.
We shall show later that MPCP-desingularizations are sufficient to
establish (\ref{duality}) for $q=1$ and arbitrary $p$ in the case of
$d$-dimensional Calabi-Yau hypersurfaces in toric varieties (see
\ref{p1}, \ref{p1cor}).
Although MPCP-desingularizations always exist, it is important
to stress that they are not sufficient
to prove (\ref{duality}) for all $p,q$, and  $d \geq 4$, because of the
following two properties which can be easily illustrated
by means of various examples:

\begin{itemize}
\item In general, a MPCP-desingularization of a
Gorenstein toroidal singularity is
not a  manifold, but a variety with Gorenstein terminal
abelian quotient singularities.

\item Cohomology and Hodge numbers of different MPCP-desingularizations
might be different.
\end{itemize}
\bigskip

It turns out that, even for $3$-dimensional Calabi-Yau manifolds,
the mirror construction inspired from the superconformal field theory demands
consideration of higher dimensional manifolds with singularities
\cite{batyrev-borisov,candelas,schimmrigk1,schimmrigk2,schimmrigk3}.
In this case, we again meet difficulties if we wish to obtain analogues
of the duality in (\ref{duality}). Let us explain them for the example
which was discussed in \cite{candelas}.

Let $E_0$ be the unique elliptic curve having an authomorphism of order
$3$ with $3$ fixed points $p_0, p_1, p_2 \in E_0$. We consider
the natural
diagonal action of $G \cong {\bf Z}/3{\bf Z}$ on $Z = E_0\times E_0 \times
E_0$. The quotient $X = Z/G$ is a singular Calabi-Yau variety whose
smooth crepant resolution $\hat{X}$ has Hodge numbers
\[ h^{1,1}(\hat{X}) = 36,\;\; h^{2,1} (\hat{X}) = 0. \]
As the mirror partner of $X$, it has been proposed the $7$-dimensional
orbifold $Y$ obtained from the quotient of the Fermat-cubic
$(W\, :\, z_0^3 + \cdots z_8^3 = 0)$ in
${\bf P}^8$ by the order 3 cyclic group
action defined by the matrix
\[ g = {\rm diag}( 1,1,1,e^{2\pi \sqrt{-1}/3},e^{2\pi \sqrt{-1}/3} ,e^{2\pi
\sqrt{-1}/3},
e^{-2\pi \sqrt{-1}/3},e^{-2\pi \sqrt{-1}/3},e^{-2\pi \sqrt{-1}/3}). \]
By standard methods, counting $G$-invariant monomials in the
Jacobian ring, one immediately verifies that $h^{4,3}(Y) = 30$.
One could expect that a crepant  resolution of singularities of $Y$
along the $3$ elliptic curves
\[ C_0 = \{ z_3 = \cdots = z_8 = 0 \} \cap Y, \]
\[ C_1 = \{ z_0 = z_1= z_2 = z_6 = z_7 = z_8 = 0 \} \cap Y, \]
\[ C_2 = \{ z_0 = \cdots = z_5 = 0 \} \cap Y \]
would give the missing  $6$ dimensions to $h^{4,3}(Y)$ in order to obtain
$36$ (this would be  the analogue of (\ref{duality})).
But also this hope must be given up  because of a very simple reason: all
singularities along $C_0, C_1, C_2$ are terminal, i.e.,
they do not
admit any smooth crepant resolution.

\begin{ques}
What could be that suitable mathematical reasoning which would give
back the missing $6$
in the above example?
\end{ques}

{}From the viewpoint of physicists, one should consider $Y$ as
an orbifold quotient of $W$ by  $G = \{ e, g, g^{-1} \}$.
By physicists' formula (\ref{phys.form}),
\[ h^{4,3}(W,G) = h^{4,3}_e(W,G) + h^{4,3}_g(W,G) +
h^{4,3}_{g^{-1}}(W,G). \]
It is clear that $h^{4,3}_g(W,G) = h^{4,3}_{g^{-1}}(W,G)$ and
$h^{4,3}_e(W,G) = h^{4,3}(Y) = 30$. So, it remains to compute
$h^{4,3}_g(W,G)$. Notice that $W^g = C_0 \cup C_1 \cup C_2$; i.e.,
$W_i(g) = C_i$ $( i =0,1,2)$. Moreover, $g$ acts on the tangent
space $T_w$ of a point $w \in W^g$ by the matrix
\[ {\rm diag}( 1,e^{2\pi \sqrt{-1}/3},e^{2\pi \sqrt{-1}/3} ,e^{2\pi
\sqrt{-1}/3},
e^{-2\pi \sqrt{-1}/3},e^{-2\pi \sqrt{-1}/3},e^{-2\pi \sqrt{-1}/3}). \]
Therefore, $F_i(g) = 3$ $({\rm for}\; i=0,1,2)$.
So $h^{4,3}_g(W,G) = \sum_{i =0}^2 h^{4-F_i(g),3-F_i(g)}(C_i) = 3$
and the required $6$ is indeed present!

\begin{ques}
Is there  a local version of the formula {\rm (\ref{phys.form})}
for the underlying space of
a quotient singularity extending that of \ref{general}?
\end{ques}

We shall answer both questions in Sections 5 and 6.

\section{$E$-polynomials of algebraic varieties}

In this section we recall some basic properties of the {\em
$E$-polynomials} of (not necessarily smooth or compact) {\em algebraic
varieties}. $E$-polynomials are defined by means of the mixed Hodge
structure (MHS) of rational
cohomology groups with compact supports \cite{dan.hov}.
As we shall see below, these polynomials obey to similar additive and
multiplicative laws as those of the {\em usual} Euler characteristic, which
enables us to compute all the Hodge numbers coming into question in a
very convenient way.

As Deligne shows in \cite{deligne}, the cohomology
groups $H^k(X, {\bf Q})$ of a (not necessarily smooth or compact) algebraic
variety $X$ carry a natural MHS. By similar methods, one can
determine a canonical MHS by considering $H^k_c(X, {\bf Q})$, i.e., the
cohomology groups {\em with compact supports}. Compared with
$H^k(X, {\bf Q})$, the MHS on $H^k_c(X, {\bf Q})$ presents some additional
technical advantages. One of them is the existence of the
following exact sequence:

\begin{prop}
Let $X$ be an algebraic variety and $Y \subset X$ a closed subvariety. Then
there is an exact sequence
\[ \ldots \rightarrow H^k_c(X \setminus Y, {\bf Q}) \rightarrow
H^k_c(X, {\bf Q}) \rightarrow H^k_c(Y, {\bf Q}) \rightarrow  \cdots  \]
consisting of $MHS$-morphisms.
\label{exact-s}
\end{prop}

\begin{dfn} {\rm Let $X$ be an algebraic variety over ${\bf C}$
which is not necessarily compact or smooth. Denote by
$h^{p,q}(H^k_c(X, {\bf C}))$ the dimension of the $(p,q)$-Hodge component of
the  $k$-th cohomology with compact supports.
We define:
\[e^{p,q}(X) := \sum_{k \geq 0} (-1)^k h^{p,q}(H_c^k(X, {\bf C})). \]
The polynomial
\[ E(X; u,v) : = \sum_{p,q} e^{p,q}(X) u^p v^q \]
is called the {\em E-polynomial} of $X$}.
\label{e-poly}
\end{dfn}

\begin{rem}
{\rm If the Hodge structure of $X$ in \ref{e-poly} is
{\em pure}, then  the coefficients $e^{p,q}(X)$ of the E-polynomial
of $X$ are related to the usual Hodge numbers
by $e^{p,q}(X) = (-1)^{p+q}h^{p,q}(X)$. In fact, the E-polynomial
(in the general case) can be regarded as a notional refinement of the
{\em virtual Poincar\'{e} polynomial} $E(X; -u,-u)$ and, of course,
of the {\em Euler   characteristic with compact supports} $e_c(X): =
E(X, -1,-1)$. It should be also mentioned, that $e_c(X) = e(X)$, i.e.,
that $e_c(X)$ is equal to the {\em usual} Euler characteristic of $X$
(cf. \cite{fulton}, pp. 141-142).  }
\end{rem}

Using Proposition \ref{exact-s},  one obtains:

\begin{prop}
Let $X$ be a disjoint union of locally closed subvarieties $X_i$ $(i \in I)$.
Then
\[E(X;u,v) = \sum_{i \in I} E(X_i;u,v). \]
\label{proper1}
\end{prop}

\begin{dfn}
{\rm Let $X$ be a disjoint union of locally closed
subvarieties $X_i$ $(i \in I)$. We shall  write $X_{i'} < X_i$, if
$X_{i'} \neq X_i$ and $X_{i'}$ is contained in the Zariski closure
$\overline{X}_i$ of $X_i$.}
\end{dfn}

\begin{prop}
For any $i_0 \in I$, one has
\[ E(X_{i_0}; u,v) = \sum_{k \geq 0} (-1)^k \sum_{X_{i_k} < \cdots <
X_{i_1} < X_{i_0} } E(\overline{X}_{i_k}; u,v). \]
\label{stra1}
\end{prop}

\noindent
{\em Proof.}  By \ref{proper1}, we get
\[ E(X_{i_0}; u,v) = E(\overline{X}_{i_0}; u,v) -
E(\overline{X}_{i_0} \setminus X_{i_0}; u,v). \]
Moreover,
\[ E(\overline{X}_{i_0} \setminus X_{i_0}; u,v) =
\sum_{X_{i_1} < X_{i_0}} E(X_{i_1}; u,v). \]
Repeating the same procedure for  $i_1 \in I$, we obtain:
\[ E(X_{i_1}; u,v) = E(\overline{X}_{i_1}; u,v) -
E(\overline{X}_{i_1} \setminus X_{i_1}; u,v),  \]
\[ E(\overline{X}_{i_1} \setminus X_{i_1}; u,v) =
\sum_{X_{i_2} < X_{i_1}} E(X_{i_2}; u,v), \;\; \; \mbox{\rm etc. } \dots \]
This leads to the claimed formula. \hfill $\Box$
\bigskip

Applying the K\"unneth formula, we get:

\begin{prop}
Let $\pi\,: \, X \rightarrow Y$ be a locally trivial fibering in Zariski
topology. Denote by $F$ the fiber over a closed point in $Y$. Then
\[ E(X;u,v) = E(Y;u,v) \cdot E(F;u,v). \]
\label{proper2}
\end{prop}

We shall use \ref{proper1} and \ref{proper2} in the following situation.
Let $\pi\,:\, {X}' \rightarrow X$ be a proper birational morphism of
algebraic varieties ${X}'$ and $X$. Let us further assume that ${X}'$ is
smooth and $X$ has a stratification by locally closed subvarieties
$X_i$ $(i \in I)$,  such that each $X_i$ is smooth and
the restriction of $\pi$ on $\pi^{-1}(X_i)$ is a locally trivial
fibering over $X_i$ in Zariski topology. Using
\ref{proper1} and \ref{proper2}, we can compute all Hodge
numbers of ${X}'$ as follows:

\begin{prop}
Let $F_i$ $(i \in I)$ denote the fiber over a closed point of $X_i$. Then
\[ E({X}'; u,v) = \sum_{i \in I} E(X_i; u,v) \cdot E(F_i;u,v). \]
\label{formula}
\end{prop}

We shall next deal with  the case in which
 $\pi\,:\, \tilde{X} \rightarrow X$ represents
  a crepant resolution of
singularities of an algebraic variety $X$ having only
Gorenstein singularities. The problem of  main interest
is to characterize  the $E$-polynomials $E(F_i; u,v)$ in terms of
singularities  of $X$ along the $X_i$'s. This problem will be
solved in the case  when  $X$ has Gorenstein toroidal or quotient
singularities.

\section{Local crepant toric resolutions}

We shall compute here the $E$-polynomials of the fibers of
crepant toric resolution mappings  of Gorenstein toric
singularities by using
their combinatorial description in terms of convex cones.
It is assumed that the reader is familiar with the theory of toric
varieties as it is presented, for instance, in the expository article
of Danilov \cite{danilov}, or in the books
of Oda \cite{oda} and Fulton \cite{fulton}.

Let $M$, $N$ be two  free abelian groups of rank $d$,
which are dual to each other, and let
$M_{\bf R}$ and  $N_{\bf R}$ be their  real scalar extensions.
The type of every $d$-dimensional
Gorenstein toroidal singularity can be described combinatorially
by a $d$-dimensional cone $\sigma = \sigma_{\Delta}
\subset N_{\bf R}$ which supports
a $(d-1)$-dimensional lattice polyhedron $\Delta \subset N_{\bf R}$
\cite{reid1}. This
lattice  polyhedron $\Delta$ can be defined as $\{ x \in \sigma
\mid \langle x, m_{\sigma} \rangle = 1 \}$ for some uniquely
determined element
$m_{\sigma} \in M$. Let
$\check{\sigma} \subset M_{\bf R}$ be dual to $\sigma$ and
set   ${\bf A}_{\sigma} := {\rm Spec}\,  {\bf C}
\lbrack \check{\sigma} \cap M \rbrack$. Then ${\bf A}_{\sigma}$ is a
$d$-dimensional affine toric variety with only Gorenstein singularities.
We denote by $p= p_{\sigma}$ the unique torus invariant closed point in
${\bf A}_{\sigma}$.

\begin{dfn}
{\rm A finite collection ${\cal T} = \{ \theta \}$ of simplices with
vertices in $\Delta \cap N$ is called a {\em triangulation}
of $\Delta$ if the
following properties are satisfied:

(i) if $\theta'$ is a face of $\theta \in {\cal T}$, then $\theta' \in
{\cal T}$;

(ii) the intersection of any two simplices $\theta_1', \theta_2' \in
{\cal T}$  is either empty, or a common face of both of them;

(iii) $\Delta = \bigcup_{\theta \in {\cal T}} \theta$.
}
\end{dfn}

Every triangulation ${\cal T}$ of $\Delta$ gives rise to a
partial crepant toric
desingularization
$\pi_{\cal T}\, : \, X_{\cal T} \rightarrow {\bf A}_{\sigma}$ of
${\bf A}_{\sigma}$, so that
$X_{\cal T}$ has at most  abelian quotient Gorenstein singularities.

\begin{dfn}
{\rm A simplex $\theta \subset \Delta \subset
\{ x \in N_{\bf R} \mid \langle x, m_{\sigma} \rangle =1  \}$
is called {\em regular} if its vertices form a part of a
${\bf Z}$-basis of $N$. }
\end{dfn}
\noindent
It is known (see, for instance,  \cite{oda}, Thm. 1.10, p.15) that $X_{\cal T}$
is smooth
if and only if all simplices in ${\cal T}$
are regular.

\begin{theo}
Assume that $\Delta$ admits a  triangulation ${\cal T}$ into regular
simplices; i.e., that the corresponding toric variety $X_{\cal T}$
in the crepant resolution
\[ \pi_{\cal T}\, : \, X_{\cal T} \rightarrow {\bf A}_{\sigma} \]
is smooth. Then $F = \pi_{\cal T}^{-1}(p)$ can be stratified
by  affine spaces.
\label{stratification}
\end{theo}

\noindent
{\em Proof.}  Let $\theta_0$ be an arbitrary  $(d-1)$-dimensional simplex
in ${\cal T}$ with vertices  $ e_1, \ldots, e_d$.
Choose  an 1-parameter multiplicative
subgroup $G_{\omega} \subset ({\bf C}^*)^d$
whose action on ${\bf A}_{\sigma}$
is defined by a weight-vector $\omega \in \sigma \cap N$, so that
$\omega = \omega_1 e_1 + \cdots + \omega_d e_d$, where
$\omega_1, \ldots, \omega_d$ are positive integers. The action of
$G_{\omega}$ on ${\bf A}_{\sigma}$  extends naturally to an action
on $X_{\cal T}$.
If $\{ \theta_0, \theta_1, \ldots, \theta_s \}$ denotes  the set of all
$(d-1)$-dimensional simplices in ${\cal T}$, then
$\sigma = \bigcup_{i =0}^s \sigma_{\theta_i}$, and
$X_{\cal T}$ is canonically covered by the corresponding
$G_{\omega}$-invariant open subsets $U_0, \ldots, U_s$, so that
$U_i \cong {\bf C}^d$.  Denote by $p_i$ $(i = 0,1, \ldots, s)$ the
unique torus invariant point in $U_i$. We assume that $\omega$ has been
already chosen in such a way, that $p_i$ is the unique $G_{\omega}$-invariant
point in $U_i$. We consider a multiplicative parameter $t$ on $G_{\omega}$
for which  the action of $G_{\omega}$ on $U_0$ is defined as follows:
\[ t \cdot (x_1, \dots, x_d): = (t^{\omega_1}x_1, \ldots, t^{\omega_d}x_d). \]
Furthermore, we set:
\[ X_i : = \{ x = (x_1, \ldots, x_d) \in X_{\cal T} \mid \lim_{t \rightarrow
\infty}
t(x) = p_i \}.  \]
Since $\pi_{\cal T}(p_i) = p$, we have $X_i \subset F$. By compactness of $F$,
for every point $x \in F$, there exists $\lim_{t \rightarrow \infty}
t(x)$ which is a $G_{\omega}$-invariant point; i.e.,
$\lim_{t \rightarrow \infty}
t(x) = p_i$ for some $i$ $(0 \leq i \leq s)$. So
$\bigcup_{i =0}^s X_i = F$. Obviously, $X_i \subset U_i$.
Moreover,  $X_i \cap X_j = \emptyset$ for $i \neq j$.

If we now choose appropriate  torus coordinates
$y_1, \ldots, y_d$  on $U_i$, so that $G_{\omega}$ acts by
\[ t\cdot (y_1, \ldots, y_k, y_{k+1}, \ldots, y_d) =
(t^{\lambda_1} y_1, \ldots, t^{\lambda_k} y_k, t^{\lambda_{k+1}}y_{k+1},
\ldots, t^{\lambda_d}y_d ) \]
with $\lambda_1, \ldots, \lambda_k$ positive and
$\lambda_{k+1}, \ldots, \lambda_d$ negative, $X_i$ is defined by
the equations $y_1 = \ldots = y_k = 0$. Therefore, $X_i$ is isomorphic to
an affine space.
\hfill $\Box$
\medskip

Let $l(k\Delta)$ denote the number of lattice points of $k\Delta$.
Then the {\em Ehrhart power series}
\[ P_{\Delta}(t): = \sum_{k \geq 0} l(k\Delta) t^k \]
can be considered as a numerical characteristic
of the toric singularity at $p_{\sigma}$.
It is well-known (see, for instance,
\cite{batyrev1}, Thm 2.11, p.356)  that $P_{\Delta}(t)$ can be
always written in the form:
\[ P_{\Delta}(t) = \frac{\psi_0(\Delta) + \psi_1(\Delta)t + \cdots
+ \psi_{d-1}(\Delta)t^{d-1}}{(1-t)^d}, \]
where $\psi_0(\Delta) = 1$
and $\psi_1(\Delta), \ldots, \psi_{d-1}(\Delta)$ are certain
nonnegative integers.

\begin{theo}
Let $\Delta$ be as in {\rm \ref{stratification}}, and
$F = \pi_{\cal T}^{-1}(p)$.
Then the cohomology groups
$ H^{2i}(F, {\bf C}), \;\; i = 0, \ldots, d-1$
are  generated by the $(i,i)$-classes of algebraic cycles, and
$H^{j}_c(F, {\bf C}) = 0$ for odd values of $j$ . Moreover,
$h^{i,i}(F) = \psi_i(\Delta) \;\; i = 0, \ldots, d-1$.
In particular, the  dimensions  $h^{i,i}(F) = {\rm dim}\, H^{2i}(F,
{\bf C})$ $( 0 \leq i \leq d-1)$
do not depend on the choice of the triangulation ${\cal T}$.
\label{invariants}
\end{theo}

\noindent
{\em Proof.} The first statement follows immediately from
Theorem \ref{stratification}.
Since $F$ is compact, we have $H^i(F,{\bf C}) = H^i_c(F, {\bf C})$.
Therefore, it is sufficient to compute  the $E$-polynomial
\[ E(F; u,v) = \sum_{p,q} e^{p,q}(F) u^p v^q. \]
Since $X_{\cal T}$ is a toric variety, it admits a natural stratification
by strata which are isomorphic to algebraic tori $T_{\theta}$ corresponding
to regular subsimplices $\theta \in {\cal T}$, such that
\[ \mbox{\rm dim}\, T_{\theta} + \mbox{\rm dim}\, \theta  = d-1. \]
The natural stratification of $X_{\cal T}$ induces a stratification of $F$.
Notice that $\pi_{\cal T} (T_{\theta}) = p_{\sigma}$ (i.e.,
$T_{\theta} \in F$)  if and only if $\theta$ does not belong to the
boundary of $\Delta$.
If $a_i$ denotes the number of $i$-dimensional regular simplices
of ${\cal T}$ which do not belong to the boundary of $\Delta$, then
$a_i$  can be identified with the  number of
$(d-1-i)$-dimensional tori in the natural stratification of
$\pi_{\cal T}^{-1}(p)$. By \ref{proper1},  we get:
\[ E(F; u,v) = \sum_{\pi_{\cal T} (T_{\theta}) = p} E(T_{\theta}; u, v). \]
Since    $E(({\bf C}^*)^k; u,v) = (uv -1 )^k$, we
obtain
\[ E(F; u,v) = a_0(u v -1 )^{d-1} +
a_1(u v -1 )^{d-2} + \cdots + a_{d-1}. \]
Now we compute $P_{\Delta}(t)$ by using the numbers  $a_i$.
If  $\theta \in {\cal T}$ is  a $i$-dimensional regular simplex, then
\[ l(k\theta) = {k+i \choose k}; \;\;\;\;\;\; {\rm i.e.,}\;  \;\;\;\;\;
P_{\theta}(t) = \frac{1}{(1-t)^{i+1}}. \]
Applying the usual inclusion-exclusion principle
for the counting of lattice points of $k\Delta$, we obtain:
\[ l(k\Delta) = \sum_{i =0}^{d-1}
\sum_{{\rm dim}\, \theta=d-1-i} (-1)^i l(k\theta), \]
where $\theta$ runs over all regular simplices in ${\cal T}$ which do not
belong to the boundary of $\Delta$. Thus,
\[ P_{\Delta}(t) =
\frac{a_{d-1}}{(1-t)^d} - \frac{a_{d-2}}{(1 - t)^{d-1}} + \cdots +
(-1)^{d-1} \frac{a_0}{(1 - t)} \]
and  the polynomial
\[ \psi_0(\Delta) + \psi_1(\Delta)t + \cdots \psi_{d-1}(\Delta)t^{d-1}
= P_{\Delta}(t) (1-t)^d \]
is equal to
\[ a_{d-1} + a_{d-2}(t -1) + \cdots + a_0(t-1)^{d-1}. \]
The latter coincides with the $E$-polynomial $E(F; u,v)$ after
making the substitution  $t = uv$. Hence,
$\psi_i(\Delta) = e^{i,i}(F)$ ($0 \leq i \leq d-1$).
\hfill $\Box$

\begin{dfn}
{\rm Let $\Delta$ be a $(d-1)$-dimensional lattice polyhedron
defining a  $d$-dimensional Gorenstein toric
singularity $p \in {\bf A}_{\sigma}$. Then
\[ S(\Delta;uv): =  \psi_0(\Delta) + \psi_1(\Delta)uv +
\cdots + \psi_{d-1}(\Delta)(uv)^{d-1} \]
will be called the {\em $S$-polynomial} of the Gorenstein
toric singularity at $p$. }
\end{dfn}

\begin{coro}
The Euler number $e(F)$ equals $S(\Delta, 1) = (d-1)! {\rm vol}(\Delta)$.
\label{eu.number}
\end{coro}

\noindent
{\em Proof.}  By definition of $P_{\Delta}(t)$,
\[ \psi_0(\Delta) + \psi_1(\Delta) + \cdots + \psi_{d-1}(\Delta)
= (d-1)! {\rm vol}(\Delta). \]
Obviously, the left hand side equals $e(F)$.
\hfill $\Box$

\begin{rem}
{\rm It is known  that the  coefficient
$\psi_{d-1}(\Delta)$ equals $l^*(\Delta)$, i.e.,  the  number
of rational points in the interior of
$\Delta$ (see \cite{dan.hov}, pp. 292-293). }
\label{lead}
\end{rem}

\section{Gorenstein
quotient singularities}

Let $G$ be a finite  subgroup of $SL(d, {\bf C})$. We shall use the fact
that any element $ g \in G$ is obviously conjugate to a diagonal matrix.

\begin{dfn}
{\rm If  an element $g \in G$ is conjugate to
\[  {\rm diag}( e^{2\pi \sqrt{-1}\alpha_1}, \ldots, e^{2\pi \sqrt{-1}\alpha_d}
)
\]
with $\alpha_i \in {\bf Q} \cap [0,1)$, then the sum
\[ wt(g): = \alpha_1 + \cdots + \alpha_d \]
will be called the {\em weight} of the element $g \in G$.
The number  $ht(g): = {\rm rk} (g - e)$ will be called the {\em height} of
$g$. }
\end{dfn}

\begin{prop}
For any $g \in G$, one has
\[ wt(g) + wt(g^{-1}) = ht(g) = ht(g^{-1}). \]
\label{dualit}
\end{prop}

\noindent
{\em Proof.}  Let
$g = {\rm diag}( e^{2\pi \sqrt{-1}\alpha_1}, \ldots, e^{2\pi \sqrt{-1}\alpha_d}
)$,
$g^{-1} = {\rm diag}( e^{2\pi \sqrt{-1}\beta_1}, \ldots, e^{2\pi
\sqrt{-1}\beta_d} )$.
Then $ht(g)$ equals the number of nonzero elements in
$\{ \alpha_1, \ldots, \alpha_d \}$. On the other hand,
$\alpha_i + \beta_i = 1$ if $\alpha_i \neq 0$, and
$\alpha_i + \beta_i = 0$ otherwise. Hence
$\sum_{i =1}^d  (\alpha_i + \beta_i) = ht(g)$.
\hfill $\Box$

\begin{conj}
{\rm (strong McKay correspondence)}
Let  $G \subset SL(d, {\bf C})$ be a finite group. Assume
that $X = {\bf C}^d/G$ admits a smooth crepant desingularization
$\pi \,: \, \hat{X}
\rightarrow X$ and $F:= \pi^{-1}(0)$.  Then $H^*(F, {\bf C})$ has a
basis consisting of classes of algebraic cycles $Z_{\{g\}} \subset F$
which are
in $1$-to-$1$ correspondence with the conjugacy classes $\{g\}$ of $G$,
so that
\[ {\rm dim}\, H^{2i}(F, {\bf C}) =
\# \{ \mbox{\rm conjugacy classes $\{g\} \subset G$,
such that $wt(g) = i$} \}. \]
\label{strong1}
\end{conj}

Now we give several evidences in support of Conjecture \ref{strong1}.

\begin{theo}
Let  $G \subset SL(d, {\bf C})$ be a finite abelian group. Suppose
that $X = {\bf C}^d/G$ admits a smooth crepant toric desingularization
$\pi \,: \, \hat{X}
\rightarrow X$ and $F: = \pi^{-1}(0)$.  Then $H^*(F, {\bf C})$ has a
basis consisting of classes of algebraic cycles $Z_g \subset F$
which are
in $1$-to-$1$ correspondence with the elements $g$ of $G$,  so that
\[ {\rm dim}\, H^{2i}(F, {\bf C}) =
\# \{ \mbox{\rm elements $g \in G$,
such that $wt(g) = i$} \}. \]
In particular, the Euler number of $F$ equals $\mid G \mid$.
\label{strong}
\end{theo}

\noindent
{\em Proof.}  Let $N \subset {\bf R}^d$ be the free abelian group
generated by ${\bf Z}^d \subset {\bf R}^d$ and all
vectors $(\alpha_1, \ldots, \alpha_d)$ where
$g = {\rm diag}( e^{2\pi \sqrt{-1}\alpha_1}, \ldots, e^{2\pi \sqrt{-1}\alpha_d}
)$ runs
over all the elements of $G$. Then  $N$ is a  full sublattice
of ${\bf R}^d = {\bf N}_{\bf R}$,  ${\bf Z}^d$ is a subgroup
of finite index in $N$, and $N/{\bf Z}^d$ is canonically isomorphic to $G$.
 Let $M = {\rm Hom}(N, {\bf Z})$. We identify
 ${\bf Z}^d$ with ${\rm Hom}({\bf Z}^d, {\bf Z})$ by using the dual
 basis. $M$ is a canonical sublattice of ${\bf Z}^d$ and
 can be identified with the set of all Laurent monomials in
variables $t_1, \ldots, t_d$ which are $G$-invariant.
Therefore, the cone $\sigma$ defining the affine toric variety
$X = {\bf A}_{\sigma}$ is the positive $d$-dimensional octant
${\bf R}_{\geq 0}^d \subset {\bf R}^d = N_{\bf R}$.
Furthermore, the element $m_{\sigma} \in M$,
which was mentioned at the beginning of the previous section,
 equals $(1,\ldots, 1) \in {\bf Z}^d$.
Now if $S: = {\bf C} \lbrack \sigma \cap N \rbrack $ and if for any
$x \in \sigma \cap N$, we define a {\em degree}
${\rm deg}\, x : = \langle x, m_{\sigma} \rangle$,
$S$ becomes a graded ring, so that
\[ n_1 := (1, 0, \ldots ,0), \ldots, n_d : = (0,\ldots, 0,1) \]
form a regular sequence of  elements of degree $1$ in $S$.
This means that  $S/(n_1, \ldots, n_d)$ has a monomial basis
corresponding to those elements of $N$ which are not in ${\bf Z}^d$.
The element $(\alpha_1, \ldots, \alpha_d) \in N$ corresponds
precisely to
the element $g = {\rm diag}( e^{2\pi \sqrt{-1}\alpha_1}, \ldots, e^{2\pi
\sqrt{-1}\alpha_d}
) \in G$. Moreover,
\[ \langle (\alpha_1, \ldots, \alpha_d), m_{\sigma} \rangle = w(g). \]
Thus, the Poincar\'{e} series of the quotient ring $S/(n_1, \ldots, n_d)$
equals
\[ \psi_0(\Delta) + \psi_1(\Delta)t + \cdots + \psi_{d-1}(\Delta)t^{d-1} \]
with coefficients
\[ \psi_i(\Delta) =   \# \{ \mbox{\rm elements $g \in G$
such that $wt(g) = i$} \} \]
and  $\sigma = \sigma_{\Delta}$ as in \S 4.
The proof is completed after making use of Theorem \ref{invariants}
and Corollary \ref{eu.number}.
\hfill $\Box$

\begin{exam}
{\rm For an abelian finite group $G \subset SL(3, {\bf C})$, the quotient
$X = {\bf C}^3 /G$ admits always smooth crepant toric desingularizations
coming from the full triangulations of the corresponding triangle $\Delta$
which is determined by $n_1, n_2, n_3$. All these triangulations contain only
regular simplices and each of them differs from another one by finitely
many elementary transformations (cf. \cite{oda}, Prop. 1.30 (ii)). In
particular, if $G$ is a cyclic group generated  by
\[ {\rm diag}( e^{\frac{2\pi \sqrt{-1}\lambda_1}{|G|} },
e^{\frac{2\pi \sqrt{-1} \lambda_2}{|G|}},
e^{\frac{2\pi \sqrt{-1}\lambda_3}{|G|} }) \]
with
\[ 0 < \lambda_1, \lambda_2, \lambda_3 < | G |, \;
\lambda_1 + \lambda_2 + \lambda_3 = | G |,\;
{\rm gcd}(\lambda_1, \lambda_2, \lambda_3) =1, \]
then:
\[ {\rm dim}\, H^0(F, {\bf C}) = 1, \;{\rm dim}\, H^1(F, {\bf C}) =
{\rm dim}\, H^3(F, {\bf C}) =  {\rm dim}\, H^5(F, {\bf C}) = 0, \]
\[ {\rm dim}\, H^2(F, {\bf C}) = \frac{1}{2} \left( |G| + \sum_{i =1}^3
{\rm gcd}(\lambda_i, |G|) \right) - 2 \]
and
\[ {\rm dim}\, H^4(F, {\bf C}) = \frac{1}{2} \left( |G| - \sum_{i =1}^3
{\rm gcd}(\lambda_i, |G|) \right) +1. \]}
\end{exam}

\begin{prop}
The Conjecture \ref{strong1} is true for $d \leq 3$.
\end{prop}

\noindent
{\em Proof.}  If $d = 2$, then $wt(g) = 1$ unless $g = e$. The
number of the conjugacy classes with weight $1$ is equal to
the number of nontrivial irreducible representations of $G$. Since
the exceptional locus $F$ of a crepant resolution is a tree of
rational curves, ${\rm dim}\, H^0(F, {\bf C})$= 1, and
${\rm dim}\, H^2(F, {\bf C})$ is the number of
irreducible components of $F$. By the classical McKay correspondence
\cite{gonzalez,knorrer,mckay},
we obtain the statement \ref{strong1}.

If $d=3$, we use the result of Roan \cite{roan3} about the existence of
crepant resolutions and the Euler number of the exceptional locus. Let
$F$ be the exceptional locus over $0$ of a crepant resolution
$\pi \, : \, \hat{X} \rightarrow X$. Then $F$ is a strong
deformation retract of
$\hat{X}$; i.e., $H^i(F, {\bf C}) = H^i(\hat{X}, {\bf C})$. On the other
hand, $H^4(\hat{X}, {\bf C})$ is Poincar\'{e} dual to $H_c^2(\hat{X}, {\bf
C})$.
Note that ${\rm dim}\, H^4(F, {\bf C})$  is nothing but
the number of irreducible
$2$-dimensional components of $F$. Since  $H^2(\hat{X}, {\bf Z})$ is
isomorphic to the Picard group of $\hat{X}$, ${\rm dim}\,
H^2(\hat{X}, {\bf C})$ is equal to the number of $\pi$-exceptional divisors.
Moreover, the subspace $H^2_c(\hat{X}, {\bf C}) \subset
H^2(\hat{X}, {\bf C})$ is spanned exactly by the classes of those exceptional
divisors whose image under $\pi$ is $0$. Therefore,
\[ {\rm dim}\, H^2(\hat{X}, {\bf C}) -
{\rm dim}\, H^4(\hat{X}, {\bf C}) = \]
\[  = \# \{ \mbox{\rm exceptional divisors $E \subset \hat{X}$,
 such that $ \pi (E) $ is a curve on $X$} \}.   \]
By the classical McKay correspondence in dimension $2$,
\[ {\rm dim}\, H^2(\hat{X}, {\bf C}) -
{\rm dim}\, H^4(\hat{X}, {\bf C}) = \]
\[ = \# \{ \mbox{\rm conjugacy classes $\{g\} \subset G$,
such that $wt(g) = 1$ and $ht(g) = 2$} \}. \]
By \cite{roan3},
\begin{equation}
 1 + {\rm dim}\, H^2(\hat{X}, {\bf C}) + {\rm dim}\, H^4(\hat{X}, {\bf C})
= \# \{ \mbox{\rm all conjugacy classes $\{g\} \subset G$}  \}.
\label{euler3}
\end{equation}
By \ref{dualit},
\[ \# \{ \mbox{\rm conjugacy classes $\{g\} \subset G$,
with  $wt(g) = 1$ and $ht(g) = 3$} \} = \]
\[ = \# \{ \mbox{\rm conjugacy classes $\{g\} \subset G$,
with  $wt(g) = 2$ and $ht(g) = 3$} \}. \]
Hence,
\[ \# \{ \mbox{\rm conjugacy classes $\{g\} \subset G$,
with  $wt(g) = 2$ and $ht(g) = 3$ } \} =
{\rm dim}\, H^4(\hat{X}, {\bf C}). \]
Notice that if $wt(g)= 2$, then the height of $g$ must be equal to $3$.
Thus,
\[ {\rm dim}\, H^4(F, {\bf C}) =
\# \{ \mbox{\rm conjugacy classes $\{g\} \subset G$,
such that  $wt(g) = 2$ } \} . \]
Finally,
\[ {\rm dim}\, H^2(F, {\bf C}) =
\# \{ \mbox{\rm conjugacy classes $\{g\} \subset G$,
such that  $wt(g) = 1$ } \}  \]
follows immediately from (\ref{euler3}).
\hfill $\Box$
\bigskip

\begin{dfn}
{\rm Let $G$ be a finite subgroup of $SL(d, {\bf C})$
and $0 \in {\bf C}^d/G$  the corresponding $d$-dimensional Gorenstein toric
singularity. If we denote by $\psi_i(G)$ the
number of the conjugacy classes of $G$ having the weight $i$,
then
\[ S(G;uv):  =  \psi_0(G) + \psi_1(G)uv +
\cdots + \psi_{d-1}(G)(uv)^{d-1} \]
will be called the {\em $S$-polynomial} of the regarded Gorenstein
quotient  singularity at $0$. }
\end{dfn}

\begin{dfn}
{\rm Let $G$ be a finite subgroup of $SL(d, {\bf C})$
and $0 \in {\bf C}^d/G$ the corresponding $d$-dimensional Gorenstein toric
singularity. If we denote by $\tilde{\psi}_i(G)$ the
number of the conjugacy classes of $G$ having the weight $i$ and
the height $d$, then
\[ \tilde{S}(G;uv): =  \tilde{\psi}_0(G) + \tilde{\psi}_1(G)uv +
\cdots + \tilde{\psi}_{d-1}(G)(uv)^{d-1} \]
will be called the {\em $\tilde{S}$-polynomial} of the Gorenstein
quotient  singularity at $0$. }
\end{dfn}

By \ref{dualit}, we easily obtain:

\begin{prop}
The $\tilde{S}$-polynomial satisfies the following reciprocity
relation:
\[ \tilde{S}(G; uv) = (uv)^d \tilde{S}(G; (uv)^{-1}). \]
\label{dualit1}
\end{prop}

\section{String-theoretic  Hodge numbers}

Let $X$ be a compact $d$-dimensional Gorenstein variety
with ${\rm Sing}\,X$ consisting of at most toroidal or quotient singularities.

\begin{dfn}
{\rm Let $x \in {\rm Sing}\,X$. We say that the $d$-dimensional
singularity at $x$ has the {\em splitting codimension $k$}, if $k$
is the maximal number for which
the analytic germ at $x$ is locally isomorphic to the
product of ${\bf C}^{d-k}$ and a $k$-dimensional
toric singularity defined by
 a $(k-1)$-dimensional lattice polyhedron $\Delta'$ or, correspondingly,
to  the product  of ${\bf C}^{d-k}$ and the underlying space
${\bf C}^k/G'$ of a
$k$-dimensional quotient singularity defined by a finite subgroup
$G' \subset SL(k, {\bf C})$. For simplicity, we also
say that the singularity at $x$ is defined by $\Delta'$, or by $G'$.}
\end{dfn}

Using standard arguments, we can easily show that $X$ is
always stratified by locally closed
subvarieties $X_i$ $(i \in I)$,  such that the germs of the singularities
of $X$ along $X_i$ are analytically isomorphic to  that of a
Gorenstein  toric singularity defined by means of
a $(k-1)$-dimensional lattice
polytope $\Delta_i$ or to that of a quotient
singularity defined by means of a finite subgroup $G_i$ of $SL(k, {\bf C})$,
respectively, where $k$ denotes the splitting codimension of singularities
on $X_i$.

\begin{dfn}
We denote by $S(X_i; uv)$ the $S$-polynomial  $S(\Delta_i; uv)$
or $S(G_i; uv)$. Analogously, $\tilde{S}(X_i; uv)$ will denote
the $\tilde{S}$-polynomial $\tilde{S}(G_i; uv)$ if $X_i$ has only
Gorenstein quotient singularities.
\end{dfn}

\begin{dfn}
{\rm Suppose that $X$ has at most quotient Gorenstein singularities.
A stratification $X = \bigcup_{i \in I} X_i$, as above, is called
{\em canonical}, if for every $i \in I$ and every $x \in X_i$, there
exists an open subset  $U \cong {\bf C}^d/G_i$ in $X$ and
an element $g \in G_i$, such that
$\overline{X_i} \cap U = ({\bf C}^d)^g/C(g)$,  where $({\bf C}^d)^g$ is the
set of $g$-invariant points of ${\bf C}^d$.}
\label{stratif}
\end{dfn}

\begin{rem}
{\rm An algebraic variety is called {\em V-variety} if it has at most
quotient singularities. A {\em Gorenstein $V$-variety} (abbreviated
{\em $GV$-variety}) is a $V$-variety having at most Gorenstein
quotient singularities. The notion
of $V$-variety (or $V$-manifold) was first
introduced by  Satake \cite{satake}. The existence and the uniqueness
of the canonical stratification for a $V$-variety was proved
by Kawasaki in \cite{kawasaki}. ( Note that our {\em canonical} stratification
in \ref{stratif} is not the first, but the second stratification of $X$
defined by Kawasaki in \cite{kawasaki}, p. 77.) }
\end{rem}

\begin{prop}
Suppose that $X$ is a $GV$-variety  and
$X = \bigcup_{i \in I} X_i$ is its  canonical stratification. Then
for any $i_0 \in I$, one has:
\[ S(X_{i_0}; uv) = \tilde{S}(X_{i_0}; uv) +
\sum_{X_{i_0} < X_{i_1}} \tilde{S}(X_{i_1}; uv). \]
\label{can.strat}
\end{prop}

\noindent
{\em Proof}. It is sufficient to prove the
corresponding local statement; i.e.,  we can assume, without loss of
generality,  that
$X_{i_0} = {\bf C}^k/G_{i_0}$. For simplicity, we set
$Y = {\bf C}^k$, $Z = X_{i_0}$. Denote by $\pi$ the
natural mapping $Y \rightarrow Z$. For  $g \in G_{i_0}$,
the image $Z(g) : = \pi(Y^g) \subset Z$
depends only on the conjugacy class of $g$.
Since $ht(g)$ equals the codimension of $Z(g)$ in $Z$, we
obtain
\[ S(X_{i_0}; uv) = \tilde{S}(X_{i_0}; uv) +
\sum_{X_{i_0} < X_{i_1}} \tilde{S}(X_{i_1}; uv). \]
\hfill $\Box$

\begin{coro}
Suppose that $X$ is a $GV$-variety  and
$X = \bigcup_{i \in I} X_i$ is its canonical stratification. Then
for any $i_0 \in I$ one has:
\[ \tilde{S}(X_{i_0}; uv) =
\sum_{k \geq 0} (-1)^k \sum_{X_{i_0} < \cdots < X_{i_k}} {S}(X_{i_k}; uv). \]
\label{stra2}
\end{coro}

\noindent
{\em Proof.}  By \ref{can.strat}, we have
\[ \tilde{S}(X_{i_0}; uv) = S(X_{i_0}; uv) - \sum_{X_{i_0} <
X_{i_1}} \tilde{S}(X_{i_1}; uv). \]
After that  we apply  \ref{can.strat} to $X_{i_1}$:
\[ \tilde{S}(X_{i_1}; uv) = S(X_{i_1}; uv) - \sum_{X_{i_1} <
X_{i_2}} \tilde{S}(X_{i_2}; uv), \;\;\; {\rm etc} \ldots  \]
The repetition of this procedure completes
the proof of the assertion. \hfill $\Box$

\begin{dfn}
{\rm Let $X$ be a stratified variety with at most Gorenstein
toroidal or quotient singularities. We shall  call the polynomial
\[ E_{\rm st}(X; u,v):=  \sum_{i \in I}
E(X_i; u,v) \cdot S(X_i; uv) \]
the {\em string-theoretic $E$-polynomial of $X$}.
Let us write $E_{\rm st}(X;u,v)$ in the following expanded form:
\[ E_{\rm st}(X; u,v) = \sum_{p,q} a_{p,q} u^p v^q. \]
The numbers  $h^{p,q}_{\rm st}(X): = (-1)^{p+q} a_{p,q}$
 will be called the {\em string-theoretic Hodge numbers of $X$}.
 Correspondingly,
\[ e_{\rm st}(X): = E_{\rm st}(X; -1,-1) = \sum_{p,q} (-1)^{p+q}
h^{p,q}_{\rm st}(X) \]
will be called the {\em string-theoretic Euler number of $X$}. }
\end{dfn}

\begin{rem}
{\rm If $X$ admits
a smooth  crepant toroidal desingularization $\pi\, : \,
\hat{X} \rightarrow X$, then,   by \ref{formula} and \ref{invariants},
the $E$-polynomial of $\hat{X}$ equals
\[ E(\hat{X}; u, v) = \sum_{i \in I}
E(X_i; u, v)\cdot E(F_i; u, v) \]
where $F_i$ denotes a the special fiber $\pi^{-1}(x)$ over
a point $x \in X_i$.}
\label{crep1}
\end{rem}

By \ref{crep1}, we obtain:

\begin{theo}
If $X$ admits a smooth  crepant toroidal desingularization
$\hat{X}$, then the string-theoretic Hodge numbers
$h^{p,q}_{\rm st}(X)$ coincide with the ordinary
Hodge numbers $h^{p,q}(\hat{X})$. In particular, the numbers
$h^{p,q}_{\rm st}(X)$ are nonnegative and satisfy the Poincar\'{e} duality
$h^{p,q}_{\rm st}(X) = h^{d-p,d-q}_{\rm st}(X)$.
\end{theo}

The next theorem will play an important  role in the forthcoming statements:

\begin{theo}
Suppose that $X$ is a $GV$-variety  and
$X = \bigcup_{i \in I} X_i$ denotes its canonical stratification. Then
\[ E_{\rm st}(X; u,v) =  \sum_{i \in I}
E(\overline{X}_i; u,v) \cdot \tilde{S}(X_i; uv). \]
\label{second.f}
\end{theo}

\noindent
{\em Proof.}  By \ref{stra1}, we get
\[ E(X_{i_0}; u,v) = \sum_{k \geq 0} (-1)^k \sum_{X_{i_k} < \cdots <
X_{i_1} < X_{i_0} } E(\overline{X}_{i_k}; u,v). \]
Therefore,
\[E_{\rm st}(X; u,v) =  \sum_{i_0 \in I}
\left( \sum_{k \geq 0} (-1)^k \sum_{X_{i_k} < \cdots <
X_{i_1} < X_{i_0} } E(\overline{X}_{i_k}; u,v) \right)
\cdot S(X_{i_0}; uv) = \]
\[  = \sum_{i_k \in I} E(\overline{X}_{i_k}; u,v) \cdot
\left( \sum_{k \geq 0} (-1)^k \sum_{X_{i_k} < \cdots <
X_{i_1} < X_{i_0} } S(X_{i_0}; uv) \right). \]
By \ref{stra2}, we have
\[ \tilde{S}(X_{i_k}; uv) = \sum_{k \geq 0} (-1)^k \sum_{X_{i_k} < \cdots <
X_{i_1} < X_{i_0} } S(X_{i_0}; uv). \]
This implies the required  formula.
\hfill $\Box$

\begin{coro}
Suppose that $X$ is a $GV$-variety.
Then the numbers
$h^{p,q}_{\rm st}(X)$ are nonnegative and satisfy the Poincar\'{e} duality
$h^{p,q}_{\rm st}(X) = h^{d-p,d-q}_{\rm st}(X)$.
\label{main.prop}
\end{coro}

\noindent
{\em Proof.}  Since $\overline{X}_i$ itself is a $V$-variety, one has
$h^{p,q}(\overline{X}_i) \geq 0$, as well as the Poincar\'{e} duality
\[ E(\overline{X}_i; u,v) = (uv)^{{\rm dim}\,\overline{X}_i}
E(\overline{X}_i; u^{-1},v^{-1}). \]
On the other hand, by \ref{dualit1}, we obtain
\[ \tilde{S}(X_i; uv) = (uv)^{{\rm dim}\,\overline{X}_i}
\tilde{S}(X_i; (uv)^{-1}). \]
This implies
\[ E_{\rm st}(X; u,v) = (uv)^{{\rm dim}\,X}
E_{\rm st}(X; u^{-1},v^{-1}). \]
Since $\tilde{S}(X_i; uv)$ is a polynomial of $uv$ with
nonnegative coefficients, we conclude that $h^{p,q}_{\rm st}(X) \geq 0$.
\hfill $\Box$

\begin{theo}
Suppose that $X$ has at most toroidal Gorenstein singularities. Let
$\pi \, : \, \hat{X} \rightarrow X$ be a $MPCP$-desingularization of
$X$. Then
\[ E_{\rm st}(X;, u,v) = E_{\rm st}(\hat{X}; u,v). \]
Moreover,
\[  h^{p,1}_{\rm st}(X) = h^{p,1}(\hat{X}), \;\;\;\; \mbox{\rm for all $p$}. \]
\label{MPCP-desing}
\end{theo}

\noindent
{\em Proof. } Let $X = \bigcup_{i \in I} X_i$ be a stratification of
$X$, such that
\[ E_{\rm st}(X; u,v) =  \sum_{i \in I}
E(X_i; u,v) \cdot S(\Delta_i; uv)  \]
and $\pi\, : \, \hat{X} \rightarrow X$ be a MPCP-desingularization of
$X$. We set $\hat{X}_i : = \pi^{-1}(X_i)$. Then $\hat{X}_i$ has the
natural stratification by products
$X_i \times ({\bf C}^*)^{{\rm codim}\, \theta}$ induced by the
triangulation
\[ \Delta_i = \bigcup_{\theta \in {\cal T}_i} \theta. \]
Thus,
\[ E_{\rm st}(\hat{X}; u,v) =  \sum_{i \in I}
\left(\sum_{\theta \in {\cal T}_i} (uv-1)^{{\rm codim}\, \theta}
E(X_i; u,v) \cdot S(\theta; uv) \right). \]
By counting lattice points in $k\Delta_i$, we obtain
\[ S(\Delta_i; uv) = \sum_{\theta \in {\cal T}_i}
(uv-1)^{{\rm codim}\, \theta}
S(\theta; uv). \]
Hence,
\[ E_{\rm st}(X;, u,v) = E_{\rm st}(\hat{X}; u,v). \]
Since $\hat{X}$ has only terminal ${\bf Q}$-factorial singularities,
for any $\theta \in {\cal T}_i$ we obtain
\[ \psi_1(\theta) = 0;\;\; \mbox{\rm  i.e.,} \;\;
S(\theta; uv) = 1 + \psi_2(\theta) (uv)^2 + \cdots . \]
Therefore, the coefficient of $u^p v$ in $E_{\rm st}(\hat{X}; u,v)$
coincides with the coefficient of $u^p v$ in the usual
$E$-polynomial  $E(\hat{X}; u,v)$. As  $\hat{X}$ is a $V$-variety,
the Hodge structure in $H^*(\hat{X}, {\bf C})$ is pure, and
\[  h^{p,1}_{\rm st}(X) = h^{p,1}(\hat{X}),\; \mbox{\rm for all $p$}. \]
\hfill $\Box$

\begin{coro}
Suppose that $X$ has at most toroidal Gorenstein singularities.
Then the numbers
$h^{p,q}_{\rm st}(X)$ are nonnegative and satisfy the Poincar\'{e} duality
$h^{p,q}_{\rm st}(X) = h^{d-p,d-q}_{\rm st}(X)$.
\end{coro}

\noindent
{\em Proof.}  By \ref{MPCP-desing}, it is sufficient to prove the statement
for a $MPCP$-desingularization $\hat{X}$ of $X$. The latter follows from
\ref{main.prop}. \hfill $\Box$

\begin{theo}
Let $X$ be a smooth compact
K\"ahler  manifold of dimension $n$ over ${\bf C}$ being equipped with
an action of a finite group $G$,  such that $X$ has a $G$-invariant
volume form. Then the orbifold Hodge numbers
$h^{p,q}(X,G)$ which were defined in the introduction
coincide with the string-theoretic Hodge numbers
$h^{p,q}_{\rm st}(X/G)$.
Moreover,
\[ e(X,G) =  e_{\rm st}(X/G). \]
\label{str-eul}
\end{theo}

\noindent
{\em Proof.}  We use the canonical stratification of $Y=  X/G$:
\[ Y  = \bigcup_{i \in I}  Y_i. \]
For every stratum $Y_i$,  there exists an element $g_i \in G$,
such that $\overline{Y}_i = X^{g_i}/C(g_i)$. We note that
\[ E(\overline{Y}_i; u,v) = \sum_{p,q} (-1)^{p+q}
{\rm dim} H^{p,q}(X^{g_i})^{C(g_i)} u^pv^q. \]
Now the equality
\[ h^{p,q}_{\rm st}(X/G) = h^{p,q}(X,G) \]
follows from Theorem \ref{second.f}.
\newline
In order to get $e(X,G) =  e_{\rm st}(X/G)$, it remains to prove
the equality
\[  e(X,G) = \sum_{p,q} (-1)^{p+q} h^{p,q}(X,G). \]
We shall make use of the notation which was introduced in \S 1.
Since $\{ g \}$ expresses  a system of representatives for
$G/C(g)$ and the number of conjugacy classes of $G$ equals
\[ \frac{1}{\mid G \mid} \sum_{g \in G} \mid C(g) \mid, \]
one can rewrite the physicists Euler number (\ref{euler.phys}) as
\[ e(X, G) = \frac{1}{\mid G \mid} \sum_{g} | C(g) | \cdot
e(X^g/C(g)) = \sum_{\{ g \}} e(X^g/C(g)),  \]
where $\{g\}$ runs over all  conjugacy classes of $G$ with $g$ representing
$\{ g \}$.
We show that
\[ \sum_{p,q} (-1)^{p+q} h^{p,q}_g(X,G) =
e(X^g/C(g)). \]
This  follows from the equalities
\[ \sum_{p,q} (-1)^{p+q} h^{p,q}_g(X,G) = \sum_{i =1}^{r_g}
\sum_{p,q} (-1)^{p+q - 2 F_i(g)} h^{p - F_i(g),q - F_i(g)}_{C(g)}(X_i(g)) = \]
\[ = \sum_{p,q} (-1)^{p+q} h^{p,q}_{C(g)}(X^g) = e(X^g/C(g)). \]
\hfill $\Box$
\bigskip

\begin{coro}
Suppose that $X/G$ has a crepant desingularization $\widehat{X/G}$
and that the strong
McKay correspondence $($Conjecture \ref{strong1}$)$ holds true
for the singularities occuring along every stratum of $X/G$. Then
\[ h^{p,q}(\widehat{X/G}) = h^{p,q}_{\rm st}(X/G). \]
\label{c.des}
\end{coro}

\begin{exam}
{\rm Let us first give a $3$-dimensional example of
an orbit space (with a simple acting group) containing both abelian and
non-abelian quotient singularities, and which was proposed by F. Hirzebruch.
We consider the Fermat quintic $X = \{ [z_1, \ldots, z_5] \in {\bf P}^4
\mid \sum_{i =1}^5 z_i^5 = 0 \}$ and let the alternating group ${\cal A}_5$
act on it coordinatewise. The group ${\cal A}_5$ has five conjugacy
classes: the trivial, one consisting of all $20$ $3$-cycles, one consisting
of the $15$ products of disjoint transpositions, and two more conjugacy
classes of $5$-cycles, each of which has $12$ elements. Note that the action
of the  elements of these last two conjugacy classes is fixed point free.
Each of the $20$ $3$-cycles fixes a plane quintic and two additional points.
Correspondingly, each of the $15$ products of disjoint traspositions fixes a
plane quintic and a projective line (without common points). As
$X/{\cal A}_5$ is a Calabi-Yau variety, the generic points of the
$1$-dimensional components of ${\rm Sing}\, X/{\cal A}_5$ are compound
Du Val points \cite{reid1}. Up to the above mentioned
$40$ additional points coming from the $3$-cycles and having isotropy
groups $\cong {\bf Z}/3{\bf Z}$, there exist $175$ more fixed points
on $X$ creating (after appropriate group identifications)
{\em dissident} points on
$X/{\cal A}_5$ (we follow here  the terminology of M. Reid).
Namely, the $25$ points of the intersection locus of the  $20$ plane quintics
(with isotropy groups $\cong {\cal A}_4$), further $125$ points
lying in the intersection locus of the $15$ plane quintics
(with isotropy groups $\cong {\cal S}_3$), as well as $15 + 10 = 25$ points
coming from the intersection of the projective lines (with isotropy groups
isomorphic to the Kleinian four-group and to
${\cal S}_3$ respectively). Using Ito's results \cite{ito1,ito2}, we can
construct global crepant desingularizations $\pi\, : \,
 \widehat{X/{\cal A}_5} \rightarrow X/{\cal A}_5$. By \ref{c.des},
 $h^{p,q}(\widehat{X/{\cal A}_5}) = h^{p,q}_{\rm st}({X/{\cal A}_5})$.
 Thus, for the
 computation of $h^{p,q}(\widehat{X/{\cal A}_5})$, we just need to choose two
 representatives, say $(123)$ and $(12)(34)$, of the two non-freely acting
 conjugacy classes. We have:

 \begin{itemize}

 \item
 $h^{p,q}({X/{\cal A}_5}) =
 h^{p,q}_{\{1\}}({X,{\cal A}_5})$ equals $\delta_{p,q}$
 ( $=$ Kronecker symbol) for
 $p + q \neq 3$, $h^{p,q}_{\{1\}}({X,{\cal A}_5}) =
 1$ for $(p,q) \in \{ (3,0), (0,3) \}$ and $h^{p,q}_{\{1\}}({X,{\cal A}_5}) =
 5$ for
 $(p,q) \in \{ (2,1), (1,2) \}$;

 \item
  $h^{p,q}_{\{(123)\}} (X, {\cal A}_5)$
 equals $2$ for $(p,q) \in \{ (1,1), (2,2) \}$, $h^{p,q}_{\{(123)\}}
 (X, {\cal A}_5) =6$  for
 $(p,q) \in \{ (2,1), (1,2) \}$,  $h^{p,q}_{\{(123)\}} (X, {\cal A}_5) =
 0$ otherwise;

\item
$h^{p,q}_{\{(12)(34)\}} (X, {\cal A}_5)$ equals $2$ for $1 \leq p,q \leq 2$
and $0$ otherwise.
\end{itemize}

Thus, we get:
\[ h^{2,1}_{\rm st}(X/{\cal A}_5) = h^{1,2}_{\rm st}(X/{\cal A}_5) = 13, \]
\[ h^{1,1}_{\rm st}(X/{\cal A}_5) = h^{2,2}_{\rm st}(X/{\cal A}_5) = 5.  \]
In particular, $e(\widehat{X/{\cal A}_5}) = e_{\rm st}(X/{\cal A}_5) =
-16$, in agreement with the calculations of physicists (cf. \cite{KS}, p. 57).
}
\end{exam}

\begin{exam}
{\rm Let $X^{(n)} : = X^n/{\cal S}_n$ be the $n$-th symmetric power
of a smooth
projective surface $X$. As it is known (see, for instance,
\cite{got2}, p.54 or \cite{hirzebruch}, p.258), $X^{(n)}$ is
endowed with a canonical crepant desingularization $X^{[n]}: =
{\rm Hilb}^n(X) \rightarrow X^{(n)}$ given by the Hilbert scheme
of finite subschemes of length $n$. In \cite{got1,got2}, G\"ottsche
computed the Poincar\'{e} polynomial of $X^{[n]}$. In particular, his
formula for the Euler number gives:
\[ \sum_{n =0}^{\infty} e(X^{[n]}) t^n =
\prod_{k =1}^{\infty}(1-t^k)^{-e(X)}. \]
Using power series comparison and the above formula, Hirzebruch and
H\"ofer  gave in \cite{hirzebruch} a formal proof of the equality
$e(X^{[n]}) = e(X^{(n)}, {\cal S}_n)$. In fact, for the proof of the
validity of {\em orbifold Euler formulae} of this kind, it is enough to
check locally that the Conjecture \ref{general} of M. Reid is true (cf.
\cite{roan3},
Lemma 1). Our results \ref{str-eul} and \ref{c.des} say more: in order to
obtain the equality $  h^{p,q}(X^{[n]}) = h^{p,q}(X^{(n)},G)$
it is sufficient to verify locally our ``strong'' McKay correspondence.
The latter has been checked by G\"ottsche in \cite{got3}. The numbers
$h^{p,q}(X^{[n]})$ can be  computed by means of the Hodge
polynomial $h(X^{[n]}; u,v) := E(X^{[n]}; -u,-v)$.

If $\Pi(n)$ denotes the set of all finite series
$(\alpha) = ( \alpha_1, \alpha_2, \ldots )$ of
nonnegative integers with $\sum_i i \alpha_i = n$, then the conjugacy
class of a permutation $\sigma \in {\cal S}_n$ is determined by
its type $(\alpha) = ( \alpha_1, \alpha_2, \ldots ) \in \Pi(n)$, where
$\alpha_i$ expresses the number of cycles of length $i$ in $\sigma$.
 G\"ottsche and Soergel
\cite{GS,got2} proved that
\[ h(X^{[n]}; u,v) = \sum_{(\alpha) \in \Pi(n)} (uv)^{n- \mid \alpha \mid}
\prod_{k =1}^{\infty} h(X^{(\alpha_k)}; u,v), \]
where $\mid \alpha \mid : = \alpha_1 + \alpha_2 + \cdots $ denotes the
sum of the members of $(\alpha) \in \Pi(n)$.

(Similar formulae can be obtained for the even-dimensional Kummer varieties
of higher order, cf. \cite{got2,got3,GS}.) }
\end{exam}

\section{Applications to quantum cohomology $\;\;\;\;\;\;\;\;\;\;\;
\;\;\; $ and mirror symmetry}

{}From now on, and throughout this section, we use the notion of
{\em reflexive polyhedron} being
introduced in \cite{batyrev1}.

\begin{prop}
Let $\Delta$ be a reflexive polyhedron of dimension $d$. Then
\[ S(\Delta,t) = (t -1)^d +
  \sum_{\begin{array}{c} {\scriptstyle 0 \leq {\rm dim}\,\theta \leq d-1}
 \\ {\scriptstyle \theta \subset \Delta} \end{array}} S(\theta,t)
 \cdot (t -1)^{{\rm dim}\, \theta^*}. \]
 \label{relation}
\end{prop}

\noindent
{\em Proof. }
Denote by $\partial \Delta$ the $(d-1)$-dimensional boundary of $\Delta$ which
is homeomorphic to $(d-1)$-dimensional sphere. Let
$l(k \cdot \partial \Delta)$ be the number of lattice points belonging
to the boundary of $k\Delta$. The reflexivity of $\Delta$ implies:
\[ l(k \cdot \partial \Delta) = \sum_{0 \leq {\rm dim}\, \theta \leq d-1}
(-1)^{d-1 - {\rm dim}\, \theta} l(k\theta), \; \; \; \mbox{\rm for
$k > 0$}.  \]
Since  the Euler number of a  $(d-1)$-dimensional sphere is $1 + (-1)^{d-1}$,
we obtain
\[ (-1)^{d-1} + (1-t)P_{\Delta}(t)\;\; = \;\;
(-1)^{d-1}  \sum_{0 \leq {\rm dim}\, \theta \leq d-1}
(-1)^{{\rm dim}\, \theta} P_{\theta}(t), \]
i.e.,
\[ (-1)^{d-1} + \frac{S(\Delta;t)}{(1-t)^d}\;\; = \;\;
(-1)^{d-1} \sum_{0 \leq {\rm dim}\, \theta \leq d-1}
(-1)^{{\rm dim}\, \theta}
\frac{S(\theta;t)}{(1-t)^{{\rm dim}\, \theta +1}}. \]
This implies the required equality.
\hfill $\Box$
\medskip

We prove  the following relation between
the polar duality of lattice polyhedra
and string-theoretic cohomology:

\begin{theo}
Let ${\bf P}_{\Delta}$ be a $d$-dimensional Gorenstein toric Fano variety
corresponding to a $d$-dimensional reflexive polyhedron $\Delta$. Then
\[ E_{\rm st}({\bf P}_{\Delta}; u, v) =
(1- uv)^{d+1} P_{\Delta^*} (uv) \]
where $\Delta^*$ is the dual reflexive polyhedron.
\end{theo}

\noindent
{\em Proof.}  ${\bf P}_{\Delta}$ has
a natural stratification being defined by the strata $T_{\theta}$,
where $\theta$ runs over all the faces of $\Delta$.
On the other hand, the Gorenstein
singularities along $T_{\theta}$ are determined
by the dual face $\theta^*$
of the dual polyhedron  $\Delta^*$ (cf. \cite{batyrev1}, 4.2.4). We set
$S(\theta^*, uv) =1$ if $\theta = \Delta$. Then
\[ E_{\rm st}({\bf P}_{\Delta}; u, v) = \sum_{\theta \subset \Delta}
E(T_{\theta}; u,v) \cdot S(\theta^*; uv). \]
Note  that
\[ E(T_{\theta}; u,v) = (uv -1)^{{\rm dim}\, \theta}, \]
and  that,  for ${\rm dim}\, \theta < d$, one has by definition:
\[ S(\theta^*; uv) = (1 - uv)^{{\rm dim}\, \theta^* +1} P_{\theta^*}(uv). \]
If we apply Proposition \ref{relation} to the
dual reflexive polyhedron $\Delta^*$, then,
using ${\rm dim}\, \theta + {\rm dim}\, \theta^* = d-1$, we obtain
the desired  formula for $E_{\rm st}({\bf P}_{\Delta}; u, v)$.
\hfill $\Box$
\bigskip

\begin{coro}
The string-theoretic Euler number of ${\bf P}_{\Delta}$ is equal to
$d!({\rm vol}\, \Delta^*)$.
\end{coro}

\begin{rem}
{\rm The quantum cohomology ring of a smooth toric variety was described
in \cite{batyrev00}. It was proved  that the usual cohomology
of a smooth toric manifold can be obtained as a limit of the quantum
cohomology ring. On the other hand, one can immediately extend
the description of the  quantum cohomology ring to arbitrary
(possibly singular) toric variety (cf. \cite{batyrev00}, 5.1).
In particular, one can easily show that
${\rm dim}\, QH^*_{\varphi}({\bf P}_{\Delta}, {\bf C})
 = d! ({\rm vol}\, \Delta^*)$,  for
any $d$-dimensional reflexive polyhedron. Comparing  dimensions, we
see that, for singular toric Fano varieties ${\bf P}_{\Delta}$,
the limit of the quantum cohomology ring is not the usual
cohomology ring, but rather the cohomology of a smooth crepant
desingularization, if  such a desingularization exists
(cf. \cite{batyrev00}, 6.5). By our general philosophy,
we should consider the string-theoretic Hodge numbers
$h^{p,p}_{\rm st}({\bf P}_{\Delta})$ as the Betti numbers
of a limit of the quantum cohomology  ring
$QH^*_{\varphi}({\bf P}_{\Delta}, {\bf C})$. }
\end{rem}
\bigskip

Let $\overline{Z}_f: = \overline{Z}_{f_1} \cap \cdots \cap \overline{Z}_{f_r}$
be a generic $(d-r)$-dimensional Calabi-Yau
complete intersection variety, which is embedded
in a Gorenstein toric Fano variety
${\bf P}_{\Delta}$ corresponding to a
$d$-dimensional reflexive polyhedron $\Delta = \Delta_1 + \cdots \Delta_r$,
where $\Delta_i$ is the Newton polyhedron of $f_i$ $(i = 1, \ldots, r)$.
Assume that the lattice polyhedra $\Delta_1, \ldots, \Delta_r$  are
defined by a {\em nef-partition} of vertices of the dual reflexive polyhedron
$\Delta^* = {\rm Conv}\{ \nabla_1, \ldots , \nabla_r\}$. (For definitions
and notations the reader is referred to \cite{batyrev-borisov,borisov}.)
Denote by $\overline{Z}_g : = \overline{Z}_{g_1}
\cap \cdots \cap \overline{Z}_{g_r}$
a generic Calabi-Yau complete  intersection variety in
the Gorenstein toric Fano variety ${\bf P}_{\nabla^*}$, which is
defined by the reflexive polyhedron $\nabla^* = {\rm Conv}\{
\Delta_1, \ldots ,
\Delta_r \}$, where $\nabla_i$ is the Newton polyhedron of
$g_i$ $(i = 1, \ldots, r)$.

\begin{conj}  {\rm (Mirror duality of string-theoretic Hodge numbers)}
The string-theoretic $E$-polynomials of $\overline{Z}_f$ and
$\overline{Z}_g$ obey to the following reciprocity law:
\[ E_{\rm st}(\overline{Z}_f; u,v) =
(-u)^{d-r}E_{\rm st}(\overline{Z}_g;u^{-1},v). \]
Equivalently, the string-theoretic Hodge numbers of $\overline{Z}_f$ and
$\overline{Z}_g$ are related to each other by:
\[ h^{p,q}_{\rm st}(\overline{Z}_f) =
h^{d-r-p,q}_{\rm st}(\overline{Z}_g),
\;\; \mbox{ {\rm for all $p,q$}}. \]
\label{symmetry}
\end{conj}

\noindent
We want to show some evidences in support  of Conjecture \ref{symmetry} for
Calabi-Yau hypersurfaces ($r =1$).

\begin{theo}
Let $\overline{Z}_f$ be a $\Delta$-regular Calabi-Yau hypersurface in
${\bf P}_{\Delta}$. Then
\[  E_{\rm st}(\overline{Z}_f; 1,v) =
\frac{S(\Delta^*;v)}{v} + (-1)^{d-1} \frac{S(\Delta;v)}{v}  + \]
\[ + \sum_{\begin{array}{c} {\scriptstyle 1 \leq {\rm dim}\,\theta \leq d-2}
 \\ {\scriptstyle \theta \subset \Delta} \end{array}}
\frac{(-1)^{{\rm dim}\, \theta-1}}{v}  \left(
 S(\theta; v)  \cdot
 S(\theta^*; v) \right) -  \]
 \[ - \sum_{\begin{array}{c} {\scriptstyle {\rm dim}\,\theta = d-1}
 \\ {\scriptstyle \theta \subset \Delta} \end{array}}
 (-1)^{d-1} \frac{S(\theta,v)}{v} -
 \sum_{\begin{array}{c} {\scriptstyle {\rm dim}\,\theta^* = d-1}
  \\ {\scriptstyle \theta^* \subset \Delta^*} \end{array}}
 \frac{S(\theta^*,v)}{v}. \]
 \label{formul0}
\end{theo}

\begin{coro}
\[ E_{\rm st}(\overline{Z}_f; 1,v) =(-1)^{d-1}
E_{\rm st}(\overline{Z}_g; 1,v).  \]
\end{coro}

At first we need the following formula:

\begin{prop} Let $\theta$ be a face of $\Delta$ and
${\rm dim}\, \theta \geq 1$. Then
\[ E(Z_{f,\theta}; 1,v) = \frac{ (v-1)^{{\rm dim}\, \theta}}{v}
+ (-1)^{{\rm dim}\, \theta-1}
\frac{S(\theta,v)}{v} . \]
\label{e-f}
\end{prop}

\noindent
{\em Proof. } It follows from the formula
of Danilov and Khovanski\^i (\cite{dan.hov}, Remark 4.6):
\[ (-1)^{{\rm dim}\, \theta -1} \sum_p e^{p,q}(Z_{f,\theta}) =
(-1)^q { n \choose q + 1 } + \psi_{q+1}(\theta). \]
\hfill $\Box$

\noindent
{\bf Proof of Theorem \ref{formul0}}. By definition,
\[ E_{\rm st}(\overline{Z}_f; 1,v) \; =  \; E(Z_{f, \Delta}; 1,v) \; + \;
\sum_{\begin{array}{c}
{\scriptstyle {\rm dim}\,\theta = d-1} \\
{\scriptstyle \theta \subset \Delta} \end{array}}
E(Z_{f, \theta}; 1,v) \;\; + \]
\[ + \sum_{\begin{array}{c}
{\scriptstyle 1 \leq {\rm dim}\,\theta \leq d-2 } \\
{\scriptstyle \theta \subset \Delta} \end{array}}
E(Z_{f, \theta}; 1,v) \cdot S(\theta^*; v). \]

Substituting the expressions which follow from  \ref{e-f}, we get:
\[   E(\overline{Z}_{f, \Delta}; 1,v) = \frac{(v-1)^{d}}{v}
+ (-1)^{d-1} \frac{S(\Delta,v)}{v} + \]
\[ +  \sum_{\begin{array}{c}
{\scriptstyle  {\rm dim}\,\theta = d-1} \\ {\scriptstyle
\theta \subset \Delta} \end{array} }
\left( (-1)^{d-2} \frac{S(\theta,v)}{v} + \frac{(v-1)^{d-1}}{v} \right) + \]
\[ +  \sum_{\begin{array}{c}
 {\scriptstyle  1 \leq {\rm dim}\,\theta \leq d-2} \\
{\scriptstyle \theta \subset \Delta} \end{array}}
\left( (-1)^{{\rm dim}\, \theta - 1} \frac{S(\theta,v)}{v} +
\frac{(v-1)^{{\rm dim}\, \theta}}{v} \right) \cdot S(\theta^*; uv). \]
It remains to use \ref{e-f} and \ref{relation}.
\hfill $\Box$
\bigskip

\begin{dfn}
{\rm For a face $\theta$ of $\Delta$, we denote by
${\bf v}(\theta)$ the normalized volume of $\theta$:
$({\rm dim}\,\theta)! {\rm vol}(\theta)$. }
\end{dfn}

\begin{coro}
Let $\Delta$ be a $d$-dimensional reflexive polyhedron. Then
\[ e_{\rm st} (\overline{Z}_f) = \sum_{i =1}^{d-2}
\sum_{{\rm dim}\, \theta = i}
(-1)^i {\bf v}(\theta)\cdot {\bf v}(\theta^*). \]
In particular,
\[ e_{\rm st} (\overline{Z}_f) = (-1)^{d-1} e_{\rm st} (\overline{Z}_g). \]
\end{coro}

We remark  that \ref{symmetry} is evident if $q =0$, because
$h^{p,0}_{\rm st} (\overline{Z}_f) = 1$,  for $q = 0, d-1$ and
$h^{p,0} _{\rm st} (\overline{Z}_f) = 0$ otherwise.
For $q = 1$ $(r=1)$,  and $p \in \{ 1,d-2\}$,
Conjecture  \ref{symmetry} is  proved by Theorem \ref{MPCP-desing}
combined with Thm. 4.4.3 from  \cite{batyrev1}.  We generalize this
for arbitrary values of $p$.

\begin{theo}
For a face $\theta$ of $\Delta$, we denote by
$l^*(\theta)$ the number of lattice points in the relative interior of
$\theta$. Assume that $d \geq 5$. Then  for $ 2 \leq  p \leq d-3$ one has
\[ h^{p,1}_{\rm st}(\overline{Z}_f) =
\sum_{{\rm codim}\, \theta = p } l^*(\theta)\cdot
l^*(\theta^*). \]
By the duality among faces, one has
\[ h^{p,1}_{\rm st}(\overline{Z}_f) =
h^{d-1-p,1}_{\rm st}(\overline{Z}_g). \]
\label{p1}
\end{theo}

\noindent
{\em Proof.}  By the Poincar\'{e} duality, it is enough  to compute
$h^{d-1-p,d-2}_{\rm st}(\overline{Z}_f) =
h^{p,1}_{\rm st}(\overline{Z}_f)$. We use
\[ E_{\rm st} (\overline{Z}_f;u,v) =
\sum_{\theta \subset \Delta} E(Z_{f,\theta};u,v)\cdot
S(\theta^*; uv). \]
By \ref{lead},
\[ S(\theta^*; uv) = l^*(\theta^*)(uv)^{{\rm dim}\, \theta^*} +
\mbox{\rm $\{$lower order  terms in $uv$$\}$}. \]
On the other hand,  by \cite{dan.hov}, Prop. 3.9,
\[e^{p,q}(Z_{f,\theta}) =0\;\; \mbox{if
$p + q > {\rm dim}\, \theta - 1 =
{\rm dim}\, Z_{f,\theta}$ and $p \neq q$}. \]
Hence,  the only possible case in which we can meet  the monomial of
type $u^{d-1-p}v^{d-2}$ within the  product
$ E(Z_{f,\theta};u,v)\cdot
S(\theta^*; u,v)$ is that occuring by consideration of
the product of the  term $l^*(\theta^*)(uv)^{{\rm dim}\, \theta^*}$ from
$S(\theta^*; uv)$ and the term
\[ e^{0, {\rm dim}\, \theta -1}(Z_{f,\theta}) v^{{\rm dim}\, \theta -1}, \]
where ${\rm dim}\, \theta^* = d - 1 - p$.
As it is known (cf. \cite{dan.hov}, Prop. 5.8.):
\[ e^{0, {\rm dim}\, \theta -1}(Z_{f,\theta})
= (-1)^{{\rm dim}\, \theta -1}l^*(\theta). \]
Therefore,
\[ h^{d-1-p,d-2}_{\rm st}(\overline{Z}_f) = l^*(\theta) \cdot
l^*(\theta^*). \]
\hfill $\Box$

\begin{coro}
Let $\hat{Z}_f$ be a MPCP-desingularization of $\overline{Z}_f$.
Assume that $d \geq 5$. Then,  for $ 2 \leq  p \leq d-3$, one has
\[ h^{p,1}(\hat{Z}_f) =
\sum_{{\rm codim}\, \theta = p } l^*(\theta)\cdot
l^*(\theta^*). \]
\label{p1cor}
\end{coro}

\noindent
{\em Proof.}  It follows from Theorem \ref{p1} and
Theorem \ref{MPCP-desing}. \hfill $\Box$

\section{Duality of string-theoretic Hodge numbers
for the Greene-Plesser construction}

In \cite{greene0,greene1} B. Greene and R. Plesser  proposed
an explicit construction of mirror pairs of Calabi-Yau orbifolds
which are obtained as abelian quotients of Fermat hypersurfaces
in weighted projective spaces. As it  was shown in  \cite{batyrev1},
5.5,  the Greene-Plesser construction can be interpreted
in terms of the polar duality of {\em reflexive simplices}.
The main purpose of this section is to
verify the mirror duality of all string-theoretic Hodge
numbers for this  construction.

{}From now on, we assume that $\Delta$ and $\Delta^*$ are
$d$-dimensional  reflexive simplices.
We shall  prove Conjecture \ref{symmetry} for
$\Delta$-regular Calabi-Yau hypersurfaces in
${\bf P}_{\Delta}$ and ${\bf P}_{\Delta^*}$.
(We remind that, for this kind of hypersurfaces and for $d = 4$,
Conjecture \ref{symmetry} was proved in
\cite{roan0,batyrev1}.)

\begin{dfn}
{\rm Let $\Theta$ be a $k$-dimensional lattice simplex. We denote
by $\tilde{S}(\Theta; uv)$ the $\tilde{S}$-polynomial of the
$(k+1)$-dimensional abelian quotient singularity defined by
$\Theta$. We denote the corresponding finite abelian subgroup
of $SL(k+1,{\bf C})$ by $G_{\Theta}$ (in the sence of \S 4,5).}
\end{dfn}

Our main statement is an immediate consequence of
the following:

\begin{theo}
Let $\overline{Z}_f$ be a $\Delta$-regular Calabi-Yau hypersurface in
${\bf P}_{\Delta}$. Then
\[  E_{\rm st}(\overline{Z}_f; u,v) =
\frac{1}{uv}\tilde{S}(\Delta^*;uv) + (-1)^{d-1} \frac{u^{d}}{v}
\tilde{S}(\Delta; u^{-1}v) + \]
\[ + \sum_{\begin{array}{c} {\scriptstyle 1 \leq {\rm dim}\,\theta \leq d-2}
 \\ {\scriptstyle \theta \subset \Delta} \end{array}}
 (-1)^{{\rm dim}\, \theta-1}  \left(
\frac{u^{{\rm dim}\, \theta}}{v} \tilde{S}(\theta; u^{-1}v)  \cdot
 \tilde{S}( \theta^*; uv) \right). \]
 \label{formul}
\end{theo}

\noindent
Indeed, if we apply Theorem \ref{formul} to the dual polyhedron
$\Delta^*$, then we get
\[  E_{\rm st}(\overline{Z}_g; u,v) =
\frac{1}{uv}\tilde{S}(\Delta;uv) + (-1)^{d-1} \frac{u^d}{v}
\tilde{S}(\Delta^*; u^{-1}v) + \]
\[ + \sum_{\begin{array}{c} {\scriptstyle 1 \leq {\rm dim}\,\theta^* \leq d-2}
 \\ {\scriptstyle \theta^* \subset \Delta^*} \end{array}}
(-1)^{{\rm dim}\, \theta^* -1}  \left(
\frac{u^{{\rm dim}\, \theta^*}}{v} \tilde{S}(\theta^*; u^{-1}v)  \cdot
 \tilde{S}( \theta; uv) \right). \]
Now the required equality
\[ E_{\rm st}(\overline{Z}_f; u,v) =(-u)^{d-1}
E_{\rm st}(\overline{Z}_g; u^{-1},v) \]
follows evidently from the $1$-to-$1$ correspondence
$\theta \leftrightarrow \theta^*$ $( 1 \leq {\rm dim}\, \theta,\,
{\rm dim}\, \theta^* \leq d-1)$  and
from the property: ${\rm dim}\, \theta + {\rm dim}\, \theta^* = d-1$.
\bigskip

For the proof of Theorem \ref{formul}, we need some preliminary
facts.

\begin{prop} Let $\theta$ be a face of $\Delta$ and
${\rm dim}\, \theta \geq 1$. Then
\[ E(Z_{f,\theta}; u,v) = \frac{ (uv-1)^{{\rm dim}\, \theta} -
(-1)^{{\rm dim}\, \theta}}{uv} + (-1)^{{\rm dim}\, \theta-1}
\left(  \sum_{\begin{array}{c}   {\scriptstyle {\rm dim}\, \tau \geq 1} \\
{\scriptstyle \tau \subset \theta} \end{array}}
\frac{u^{{\rm dim}\, \tau}}{v} \tilde{S}(\tau; u^{-1}v) \right) . \]
\label{e-ff}
\end{prop}

\noindent
{\em Proof. } By \cite{dan.hov}, Prop. 3.9, the natural mapping
\[ H^i_c(Z_{f, \theta}) \rightarrow H^{i+1}_c(T_{\theta}) \]
is an isomorphism if $i > {\rm dim}\, \theta -1$ and
surjective if  $i = {\rm dim}\, \theta -1$. Moreover,
$H^i_c(Z_{f, \theta}) = 0$ if $i < {\rm dim}\, \theta -1$.
In order to compute the mixed Hodge structure in
$H^{{\rm dim}\, \theta -1}_c(Z_{f, \theta})$, we use
the explicit description of the weight filtration in
$H^{{\rm dim}\, \theta -1}_c(Z_{f, \theta})$ (see \cite{batyrev0}).
Note that if we  choose a $\theta$-regular Laurent polynomial $f$
containing only ${\rm dim}\, \theta + 1$ monomials
associated with  vertices of
$\theta$ (such a polynomial $f$ defines a Fermat-type
hypersurface $\overline{Z}_f$ in ${\bf P}_{\theta}$), then
the corresponding Jacobian ring $R_f$ has a  monomial
basis. Thus, the weight filtration on $R_f$ can be described
in terms of the partition of
monomials in $R_f$ which is defined by the faces $\tau \subset \theta$.
To get the claimed formula, it suffices to
identify the partition of monomials in $R_f$ with
the height-partition of elements of the
finite abelian group $G_{\theta} \subset SL({\rm dim}\, \theta +1,
{\bf C})$ and its subroups $G_{\tau} \subset G_{\theta}$.

Another way to obtain the same result is to use the formulae
of Danilov and Khovanski\^i (cf. \cite{dan.hov}, \S 5.6,5.7)
which are  valid for
an arbitrary simple polyhedron $\Delta$.
\hfill $\Box$

\begin{prop}
Let $\theta$ be a face of $\Delta$ and
${\rm dim}\, \theta \geq 1$. Then
\[ S(\theta; t) = 1 + \sum_{\begin{array}{c}
{\scriptstyle {\rm dim}\, \eta \geq 1} \\ {\scriptstyle
\eta \subset \theta} \end{array} }
 \tilde{S}( \eta; t). \]
 \label{can}
\end{prop}

\noindent
{\em Proof. } It is similar to that of \ref{can.strat}. \hfill
$\Box$

\begin{prop}
We fix a face $\tau \subset \Delta$ and a face $\eta \subset \Delta^*$,
such that: $\tau$ is a face of $\eta^*$. Then
\[ \sum_{\theta,\; \tau \subset \theta \subset \eta^*}
(-1)^{{\rm dim}\, \theta} = (-1)^{{\rm dim}\, \tau}\;\;
\mbox{if $ \tau = \eta^*$ } \]
and
\[ \sum_{\theta,\; \tau \subset \theta \subset \eta^*}
(-1)^{{\rm dim}\, \theta} = 0\;\;
\mbox{if $ \tau \neq \eta^*$. } \]
\label{sum1}
\end{prop}

\noindent
{\em Proof.}
If $\eta^* = \tau$, this is  obvious.
For ${\rm dim}\, \eta^* > {\rm dim}\, \tau $,
the number of faces $\theta \subset \Delta$,  for which
$\tau \subset \theta \subset \eta^*$,  is equal to
${ {\rm dim}\, \eta^* - {\rm dim}\, \tau
\choose {\rm dim}\, \theta - {\rm dim}\, \tau  }$.
It remains to use the equality
\[  \sum_{\theta, \tau \subset \theta \eta^*}
(-1)^{ {\rm dim}\, \theta } = (-1)^{{\rm dim}\, \tau }
\left( \sum_{ i =0}^{{\rm dim}\, \eta^* - {\rm dim}\, \tau}
(-1)^i { {\rm dim}\, \eta^* - {\rm dim}\, \tau \choose i} \right) =  0. \]
\hfill $\Box$

\begin{prop}
\[  \frac{1}{uv} \tilde{S}(\Delta; uv) =
\frac{(uv)^d -1}{uv -1}\;\;\; +
\sum_{\begin{array}{c}  {\scriptstyle 1 \leq  {\rm dim}\,\tau \leq d-2} \\
{\scriptstyle \tau \subset  \Delta}  \end{array}}
\left( \frac{(uv)^{{\rm dim}\, \tau^*}  -1}{uv -1} \right) \cdot
\tilde{S}( \eta; uv). \]
\label{tilde-s}
\end{prop}

\noindent
{\em Proof.}
By Proposition \ref{relation}, we have
\[ (-1)^{d-1} + \frac{S(\Delta;t)}{(1-t)^d}\;\; = \;\;
(-1)^{d-1} \sum_{0 \leq {\rm dim}\, \theta \leq d-1}
(-1)^{{\rm dim}\, \theta}
\frac{S(\theta;t)}{(1-t)^{{\rm dim}\, \theta +1}}. \]
Applying  Proposition \ref{can} to both sides of this equality, we get
\[ (-1)^{d-1}\;\; + \;\; \frac{1}{(1-t)^d}\;\; +
\sum_{\begin{array}{c} {\scriptstyle {\rm dim}\, \tau \geq 1}
\\ {\scriptstyle \tau \subset \Delta} \end{array}}
\frac{\tilde{S}(\tau; t)}{(1-t)^d} \;\; = \]
\[ = \;\; (-1)^{d-1} \sum_{0 \leq {\rm dim}\, \theta \leq d-1}
\frac{(-1)^{{\rm dim}\, \theta}}{( 1- t)^{{\rm dim}\, \theta + 1}} \;\; + \]
\[ + \;\;  (-1)^{d-1} \sum_{0 \leq {\rm dim}\, \theta \leq d-1}
(-1)^{{\rm dim}\, \theta}
\sum_{\begin{array}{c} {\scriptstyle {\rm dim}\, \tau \geq 1}
\\ {\scriptstyle \tau \subset \theta} \end{array}}
\frac{\tilde{S}(\tau; t)}{(1-t)^{{\rm dim}\, \theta +1}}. \]
As the number of $k$-dimensional faces of $\Delta$ equals ${ d+1 \choose
k+1 }$, we have
\[ - \;\; (-1)^{d-1} \; - \;  \frac{1}{(1-t)^d} +
(-1)^{d-1} \sum_{0 \leq \theta \leq d-1}
\frac{(-1)^{{\rm dim}\, \theta}}{( 1- t)^{{\rm dim}\, \theta + 1}} \; = \]
\[ -\;\;  (-1)^{d-1} \; - \; \frac{1}{(1-t)^d} \;\; + \;\; \sum_{k =0}^{d-1}
\frac{(-1)^k}{(1-t)^{k+1}} { d+1 \choose k+1 } \; = \;
(-1)^d \frac{t^{d+1} - t}{(t-1)^{d+1}} \]
and we can deduce that:
\[ \frac{\tilde{S}(\Delta, t)}{(1-t)^d} \; +
\sum_{{\rm dim}\, \tau = d-1} \frac{\tilde{S}(\tau, t)}{(1-t)^d} \; +
\sum_{1 \leq {\rm dim}\, \tau \leq d-2}
\frac{\tilde{S}(\tau, t)}{(1-t)^d} \; = \]
\[ = (-1)^d \frac{t^{d+1} - t}{(t-1)^{d+1}}\;\;  +
\sum_{{\rm dim}\, \tau = d-1} \frac{\tilde{S}(\tau, t)}{(1-t)^d} \;\; + \]
\[ + \sum_{{\rm dim}\, \theta = d-1} \sum_{\begin{array}{c}
{\scriptstyle 1 \leq {\rm dim}\,\tau \leq d-2} \\
{\scriptstyle \tau \subset \theta} \end{array}}
\frac{\tilde{S}(\tau, t)}{(1-t)^d} \;\; + \]
\[ + \;\;  (-1)^{d-1}
\sum_{ 1 \leq {\rm dim}\, \theta \leq d-2} (-1)^{{\rm dim}\, \theta}
\sum_{\begin{array}{c} {\scriptstyle {\rm dim}\, \tau \geq 1}
\\ {\scriptstyle\tau \subset \theta} \end{array} }
\frac{\tilde{S}(\tau, t)}{(1- t)^{{\rm dim}\, \theta + 1}}. \]
The terms containing  $\tilde{S}(\tau, t)$, with
${\rm dim}\, \tau = d-1$,   have the same contribution
to  the right  and left hand sides. The coefficient of
$\tilde{S}(\tau, t)$
$( 1 \leq {\rm dim}\, \tau \leq d-2 )$ in the right hand side of the
last equality is
\[ (-1)^{d-1} \sum_{\begin{array}{c}
{\scriptstyle {\rm dim}\, \theta \leq d-2}
\\ {\scriptstyle \tau \subset \theta} \end{array}}
(-1)^{{\rm dim}\, \theta} \frac{1}{(1-t)^{{\rm dim}\, \theta + 1}} \;\; = \]
\[ = \;\; \frac{(-1)^d}{(t-1)^{d+1}} \left( t^{d- {\rm dim}\, \tau} - 1 -
(d - {\rm dim}\,\tau ) (t-1) \right). \]
Correspondingly, the coefficient of  $\tilde{S}(\tau, t)$
$(1 \leq {\rm dim}\, \tau \leq d-2)$ in the left hand side equals
\[ \frac{d - 1 - {\rm dim}\, \tau}{(1-t)^d}. \]
Finally, using ${\rm dim}\, \tau + {\rm dim}\, \tau^* = d -1$, we
obtain:
\[ \frac{\tilde{S}(\Delta, t)}{(1-t)^d} \; = \;
(-1)^d \frac{t^{d+1} - t}{(t-1)^{d+1}} \; + \;
(-1)^d \sum_{1 \leq  \tau \leq d-2}
\tilde{S}(\tau, t)
\frac{(t^{{\rm dim}\, \tau^* + 1} - t)}{ (t - 1)^{d + 1} }. \]
\hfill $\Box$

\noindent
{\bf Proof of Theorem \ref{formul}}. By definition,
\[ E_{\rm st}(\overline{Z}_f; u,v) \; =  \; E(Z_{f, \Delta}; u,v) \; + \;
\sum_{\begin{array}{c}
{\scriptstyle {\rm dim}\,\theta = d-1} \\
{\scriptstyle \theta \subset \Delta} \end{array}}
E(Z_{f, \theta}; u,v) \;\; + \]
\[ + \sum_{\begin{array}{c}
{\scriptstyle 1 \leq {\rm dim}\,\theta \leq d-2 } \\
{\scriptstyle \theta \subset \Delta} \end{array}}
E(Z_{f, \theta}; u,v) \cdot S(\theta^*; uv). \]

Substituting the expressions which were found out in \ref{e-ff} for
the $E$-polynomials of the above three summands, we get:
\[   E(Z_{f, \Delta}; u,v) = \frac{ (uv-1)^{d} - (-1)^{d}}{uv}
+ (-1)^{d-1}
\left(  \sum_{\begin{array}{c} {\scriptstyle {\rm dim}\, \tau \geq 1} \\
{\scriptstyle \tau \subset \Delta} \end{array}}
\frac{u^{{\rm dim}\, \tau}}{v} \tilde{S}(\tau; u^{-1}v) \right) ,  \]
\newline
\[ \sum_{\begin{array}{c}
{\scriptstyle  {\rm dim}\,\theta = d-1} \\ {\scriptstyle
\theta \subset \Delta} \end{array} }
E(Z_{f, \theta}; u,v) \;\; = \;\; \sum_{\begin{array}{c}
{\scriptstyle  {\rm dim}\,\theta =  d-1} \\
{\scriptstyle \theta \subset \Delta} \end{array}}
 \frac{ (uv-1)^{{\rm dim}\, \theta} -
(-1)^{{\rm dim}\, \theta}}{uv} \;\; + \]
\[ +   \sum_{\begin{array}{c}
{\scriptstyle  {\rm dim}\,\theta =  d-1} \\
{\scriptstyle \theta \subset \Delta} \end{array}} (-1)^{{\rm dim}\, \theta-1}
\left(  \sum_{\begin{array}{c}
{\scriptstyle {\rm dim}\, \tau \geq 1} \\
{\scriptstyle \tau \subset \theta} \end{array}}
\frac{u^{{\rm dim}\, \tau}}{v} \tilde{S}(\tau; u^{-1}v)  \right),  \]
and
\[  \sum_{\begin{array}{c}
 {\scriptstyle  1 \leq {\rm dim}\,\theta \leq d-2} \\
{\scriptstyle \theta \subset \Delta} \end{array}}
E(Z_{f, \theta}; u,v) \cdot S(\theta^*; uv) \;\;  =
\;\; \sum_{\begin{array}{c}
{\scriptstyle  1 \leq {\rm dim}\,\theta \leq d-2} \\
{\scriptstyle \theta \subset \Delta} \end{array}}
 \frac{ (uv-1)^{{\rm dim}\, \theta} -
(-1)^{{\rm dim}\, \theta}}{uv} \;\; + \]
\[ + \sum_{\begin{array}{c}
{\scriptstyle  1 \leq {\rm dim}\,\theta \leq d-2} \\
{\scriptstyle \theta \subset \Delta} \end{array}}  (-1)^{{\rm dim}\, \theta-1}
\left(  \sum_{\begin{array}{c}
{\scriptstyle  {\rm dim}\, \tau \geq 1} \\
{\scriptstyle \tau \subset \theta}
 \end{array}}
\frac{u^{{\rm dim}\, \tau}}{v} \tilde{S}(\tau; u^{-1}v)
 \right) \cdot
 \left( 1 + \sum_{
 \begin{array}{c} {\scriptstyle {\rm dim}\, \eta \geq 1} \\
 {\scriptstyle \eta \subset \theta^*}   \end{array}}
 \tilde{S}( \eta; uv)  \right). \]
Hence, $ E_{\rm st}(\overline{Z}_f; u,v)$
can be written as the sum  of the following $4$ terms $E_i$ $( i =1,2,3,4)$:
\newline
\[ E_1 = \sum_{\begin{array}{c}
{\scriptstyle  1 \leq {\rm dim}\, \theta} \\
{\scriptstyle \theta \subset \Delta} \end{array}}
\frac{ (uv-1)^{{\rm dim}\, \theta} -
(-1)^{{\rm dim}\, \theta}}{uv}, \]
\newline
\[ E_2 = \sum_{\begin{array}{c}
{\scriptstyle 1 \leq {\rm dim}\, \theta} \\
{\scriptstyle \theta \subset \Delta}
 \end{array}}
(-1)^{{\rm dim}\, \theta-1}
\left(  \sum_{\begin{array}{c}
{\scriptstyle  {\rm dim}\, \tau \geq 1} \\
{\scriptstyle \tau \subset \theta}
\end{array}}
\frac{u^{{\rm dim}\, \tau}}{v} \tilde{S}(\tau; u^{-1}v) \right) ,\]

\[ E_3 = \sum_{\begin{array}{c}
{\scriptstyle  1 \leq {\rm dim}\,\theta \leq d-2} \\
{\scriptstyle \theta \subset \Delta} \end{array}}
\left( \frac{ (uv-1)^{{\rm dim}\, \theta} -
(-1)^{{\rm dim}\, \theta}}{uv} \right) \cdot
\left(  \sum_{\begin{array}{c}
{\scriptstyle {\rm dim}\, \eta \geq 1} \\
{\scriptstyle \eta \subset \theta^*} \end{array}}
 \tilde{S}( \eta; uv)  \right), \]
and

\[ E_4 = \sum_{\begin{array}{c}
{\scriptstyle 1 \leq {\rm dim}\,\theta \leq d-2} \\
{\scriptstyle \theta \subset \Delta}
 \end{array}}
  (-1)^{{\rm dim}\, \theta-1}
\left(  \sum_{\begin{array}{c}
{\scriptstyle  {\rm dim}\, \tau \geq 1}  \\
{\scriptstyle \tau \subset \theta} \end{array}}
\frac{u^{{\rm dim}\, \tau}}{v} \tilde{S}(\tau; u^{-1}v) \right)
 \cdot
\left( \sum_{\begin{array}{c}
{\scriptstyle {\rm dim}\, \eta \geq 1} \\
 {\scriptstyle \eta \subset \theta^*}\end{array}}
 \tilde{S}( \eta; uv) \right). \]

By \ref{sum1},  we can simplify the multiple summation into a single
sum:
\[ E_4 = \sum_{\begin{array}{c}
{\scriptstyle  1 \leq {\rm dim}\,\theta \leq d-2} \\
{\scriptstyle \theta \subset \Delta} \end{array}}
(-1)^{{\rm dim}\, \theta-1}  \left(
\frac{u^{{\rm dim}\, \theta}}{v} \tilde{S}(\theta; u^{-1}v)  \cdot
 \tilde{S}( \theta^*; uv) \right). \]

If we make use of the combinatorial identity
 \[ \sum_{\begin{array}{c}
{\scriptstyle  1 \leq {\rm dim}\, \theta} \\
{\scriptstyle  \theta \subset \Delta} \end{array}}
 a^{{\rm dim}\, \theta} = \sum_{k =2}^{d+1} { d+1 \choose k } a^{k-1} =
 a^{-1} \left( (a+1)^{d+1} - 1 - (d+1) a \right), \]
we obtain:
 \[ E_1 = \sum_{\begin{array}{c}
{\scriptstyle   1 \leq {\rm dim}\, \theta} \\
 {\scriptstyle \theta \subset \Delta}\end{array}}
\frac{ (uv-1)^{{\rm dim}\, \theta} -
(-1)^{{\rm dim}\, \theta}}{uv} = \]
\[ = [uv(uv - 1)]^{-1} \left(  (uv)^{d+1} - 1 - (d+1)(uv -1) \right) +
d (uv)^{-1} = \frac{(uv)^d -1}{uv -1}. \]

By \ref{sum1}, we get
 \[ E_2 = \sum_{\begin{array}{c}
{\scriptstyle   1 \leq {\rm dim}\, \theta} \\
 {\scriptstyle \theta \subset \Delta} \end{array}}
(-1)^{{\rm dim}\, \theta-1}
  \sum_{\begin{array}{c}
 {\scriptstyle {\rm dim}\, \tau \geq 1} \\
{\scriptstyle \tau \subset \theta} \end{array}}
\frac{u^{{\rm dim}\, \tau}}{v} \tilde{S}(\tau; u^{-1}v)  =
(-1)^{d-1} \frac{u^d}{v} \tilde{S}(\Delta; u^{-1}v).  \]

\noindent
It remains to compute $E_3$.  As above for $E_1$, we have
 \[ \sum_{\begin{array}{c}
{\scriptstyle  1 \leq {\rm dim}\, \theta} \\
{\scriptstyle  \theta \subset \eta^*} \end{array}}
\frac{ (uv-1)^{{\rm dim}\, \theta} -
(-1)^{{\rm dim}\, \theta}}{uv} =
\frac{(uv)^{{\rm dim}\, \eta^*}  -1}{uv -1}. \]
Hence, by \ref{tilde-s},
\[ E_3 = \sum_{\begin{array}{c}
{\scriptstyle  1 \leq  {\rm dim}\,\eta \leq d-2} \\
{\scriptstyle \eta \subset  \Delta^*}
 \end{array}}
\left( \frac{(uv)^{{\rm dim}\, \eta^*}  -1}{uv -1} \right) \cdot
\tilde{S}( \eta; uv) = \frac{1}{uv}\tilde{S}(\Delta^*;uv) -
\frac{(uv)^d -1}{uv -1}. \]

Finally, we get altogether
\[  E_{\rm st}(\overline{Z}_f; u,v) = \frac{1}{uv}\tilde{S}(\Delta^*;uv)
 + (-1)^{d-1} \frac{u^d}{v}
\tilde{S}(\Delta; u^{-1}v) + \]
\[ + \sum_{\begin{array}{c}
{\scriptstyle  1 \leq {\rm dim}\,\theta \leq d-2} \\
{\scriptstyle \theta \subset \Delta}
 \end{array}}
(-1)^{{\rm dim}\, \theta-1}  \left(
\frac{u^{{\rm dim}\, \theta}}{v} \tilde{S}(\theta; u^{-1}v)  \cdot
 \tilde{S}( \theta^*; uv) \right). \]
\hfill $\Box$

\begin{exam}
{\rm The polar duality between reflexive simplices shows
(cf. \cite{batyrev1}, Thm. 5.1.1.) that the family of all smooth
Calabi-Yau hypersurfaces $X_{d+1}$ of degree $d+1$ in ${\bf P}^d$
has as its mirror partner the one-parameter family
$\{ Q_{d+1}(\lambda)/G_{d+1} \}$, where
\[  Q_{d+1}(\lambda) := \{ [z_0, \ldots, z_d ] \in {\bf P}^d \mid
\sum_{i=0}^d z_i^{d+1} - (d+1)\lambda\prod_{ i=0}^d z_i = 0 \} \]
denotes the so called {\em Dwork pencil} and
$G_{d+1}$ the acting  finite abelian group
\[ G_{d+1} := \{ (\alpha_0, \ldots, \alpha_d) \in
({\bf Z}/(d+1){\bf Z})^{d+1} \mid \prod_{i =0}^d \alpha_i = 1 \} /
\{\rm scalars \}, \]
which is abstractly isomorphic to $({\bf Z}/(d+1){\bf Z})^{d-1}$.
The moduli space ${\bf P}^1 \setminus \{ 0,1, \infty \}$ of
$\{ Q_{d+1}(\lambda)/G_{d+1} \}_{\lambda}$ can be described by
means of the parameter $\lambda^{d+1}$ (cf. \cite{greene},
\S 3.1, \cite{morrison}, \S 5, and \cite{morrison1} \S 11).

Since Conjecture \ref{symmetry} is true for the case being
under consideration, the quotient $Q_{d+1}(\lambda)/G_{d+1}$ has
the following string-theoretic Hodge numbers:
\[ h^{p,q}_{\rm st} (Q_{d+1}(\lambda)/G_{d+1}) =
h^{p,q}(Q_{d+1}(\lambda),G_{d+1}) = h^{d-1-p,q}(X_{d+1}) = \delta_{d-1-p,q},
\;\;\; \mbox{\rm for $p \neq q$}; \]
\[  h^{p,p}_{\rm st} (Q_{d+1}(\lambda)/G_{d+1}) =
h^{p,p}(Q_{d+1}(\lambda),G_{d+1}) = h^{d-1-p,p}(X_{d+1}) \;= \]
\[ = \; \sum_{i =0}^p (-1)^i { d+1 \choose i }
{ (p+1 -i)d + p \choose d } + \delta_{2p,d-1} . \]
In particular, the string-theoretic Euler number is given by:
\[ e_{\rm st} (Q_{d+1}(\lambda)/G_{d+1}) =
e(Q_{d+1}(\lambda),G_{d+1}) = - e(X_{d+1})\;  = \]
\[ = \; \frac{1}{d+1} \left(
(-1)^{d+2} \cdot d^{d+1} + 1 \right) - d - 1. \]
The first two equalities follow from Lefschetz hyperplane section
theorem and from the ``four-term formula'' (cf. \cite{hirzebruch1},
\S 2.2 ). The third one can be obtained directly by computing the
$(d-1)$-th Chern class of $X_{d+1}$. }
\end{exam}

\end{document}